\documentclass[12pt,dvips]{amsart}
\usepackage{amsfonts, amssymb, latexsym, epsfig,epic,xspace}

\newcommand\B{\reals^3_{\sum=0}}
\newcommand\BZ{\integers^3_{\sum=0}}

\newcommand\lie[1]{{\mathfrak #1}}

\newcommand\supp{{\rm supp\,}}
\newcommand\iso{{\cong}}
\newcommand\tensor{{\otimes}}
\newcommand\Tensor{{\bigotimes}}

\newtheorem{Theorem}{Theorem} 
\newtheorem{Proposition}{Proposition} 
\newtheorem{Lemma}{Lemma}
\newtheorem{Sublemma}{Sublemma}
\newtheorem*{Corollary}{Corollary}
 
\newtheorem*{Theorem*}{Theorem}
\theoremstyle{remark}
\newtheorem{Example}{Example}

\newcommand\onto{\twoheadrightarrow}
\newcommand\into{\operatorname*{\hookrightarrow}}

\newcommand\semiinfinite{semiinfinite\xspace}

\newcommand\Union{\bigcup}

\newcommand\reals{{\mathbb R}}
\newcommand\complexes{{\mathbb C}}
\newcommand\integers{{\mathbb Z}}

\newcommand\naturals{{\mathbb N}}

\newcommand\FlagsCn{{Fl(\complexes^n)}}

\theoremstyle{plain}

\newtheorem*{Conjecture}{Conjecture}

\newcommand\dfn{\bf} 

\newcommand\hive{{\tt hive}}
\newcommand\HIVE{{\tt HIVE}}
\newcommand\HONEY{{\tt HONEY}}
\newcommand\BDRY{{\tt BDRY}}
\newcommand\GLn{{{GL_n(\complexes)}}}
\newcommand\SLn{{{SL_n(\complexes)}}}

\begin{document}
\pagestyle{plain}

\title{The honeycomb model of $\GLn$ tensor products I:\\
Proof of the saturation conjecture}
\author{Allen Knutson}
\thanks{Supported by an NSF Postdoctoral Fellowship.}
\email{allenk@alumni.caltech.edu}
\address{Mathematics Department\\ Brandeis University\\ Waltham, Massachusetts}
\author{Terence Tao}
\thanks{Partially supported by NSF grant DMS-9706764.}
\email{tao@math.ucla.edu}
\address{Mathematics Department\\ UCLA\\ Los Angeles, California}
\date{\today}

\maketitle

\begin{abstract}
Recently Klyachko \cite{Kly} has given linear inequalities on triples
$(\lambda,\mu,\nu)$ of dominant weights of $\GLn$ necessary for the
the corresponding Littlewood-Richardson coefficient
$\dim (V_\lambda \tensor V_\mu \tensor V_\nu)^{\GLn}$ to be positive.
We show that these conditions
are also sufficient, which was known as the saturation conjecture.
In particular this proves Horn's conjecture from 1962, giving a recursive
system of inequalities \cite{Horn}.

Our principal tool is a new model of the Berenstein-Zelevinsky cone
for computing Littlewood-Richardson coefficients
\cite{BZ,Z}, the {\em honeycomb} model.
The saturation conjecture is a corollary of our main result, which 
is the existence of a particularly well-behaved honeycomb
associated to regular triples $(\lambda,\mu,\nu)$.
\end{abstract}


\section{The saturation conjecture}

A very old and fundamental question about the representation theory
of $\GLn$ is the following:

\begin{quotation}
 For which triples of dominant weights $\lambda,\mu,\nu$ does the
 tensor product $V_\lambda \tensor V_\mu \tensor V_\nu$ of the
 irreducible representations with those high weights contain a
 $\GLn$-invariant vector? 
\end{quotation}

Another standard, if less symmetric, formulation of the problem above 
replaces $V_\nu$ with its dual, and
asks for which $\nu$ is $V_\nu^*$ 
a constituent of $V_\lambda \tensor V_\mu$.
In this formulation one can without essential loss of generality
restrict to the case that $\lambda$, $\mu$, and $\nu^*$ are 
polynomial representations,
and rephrase the question in the language of Littlewood-Richardson
coefficients; it asks for which triple of partitions $\lambda,\mu,\nu^*$
is the Littlewood-Richardson coefficient $c^{\nu^*}_{\lambda \mu}$ positive.

It is not hard to prove (as we will later in this introduction)
that the set of such triples $(\lambda,\mu,\nu)$
is closed under addition, so forms a monoid. In this paper we prove
that this monoid is {\em saturated}, i.e. that for each triple
of dominant weights $(\lambda,\mu,\nu)$,
$$
(V_{N\lambda} \tensor V_{N\mu} \tensor V_{N\nu})^\GLn > 0 
\quad \hbox{for some $N>0$}
\qquad \implies \qquad
(V_{\lambda} \tensor V_{\mu} \tensor V_{\nu})^\GLn > 0.$$

This is of particular interest because Klyachko has recently given an answer%
\footnote{%
Klyachko gives a finite set of inequalities, that {\em as a set} are 
necessary and sufficient for this asymptotic result.
However, Chris Woodward has informed us that 
contrary to Klyachko's unproven claim in \cite{Kly}, 
the inequalities are not independent -- not all of them determine
facets of the cone. This will be the subject of inquiry of
our second paper \cite{Hon2}.}
to the general question above, which in one direction was only asymptotic
\cite{Kly}:

\begin{quotation}
 If $V_{\lambda} \tensor V_{\mu} \tensor V_{\nu}$ 
 has a $\GLn$-invariant vector, then
 $\lambda,\mu,\nu$ satisfy a certain system of linear inequalities
 derived from Schubert calculus (plus the evident linear equality
 that $\lambda+\mu+\nu$ be in the root lattice;
 in the L-R context this asks that the number of boxes in the
 partition $\nu^*$ is the number of boxes in $\lambda$ and $\mu$ together).
 Conversely, if $\lambda,\mu,\nu$ satisfy these inequalities, 
 then there exists an integer $N$ such that
 the tensor product $V_{N\lambda} \tensor V_{N\mu} \tensor V_{N\nu}$ 
 has a $\GLn$-invariant vector. 
\end{quotation}

Our saturation result completes this converse, saying
that Klyachko's inequalities completely characterize the monoid.
The survey papers \cite{F,Z} point out another important consequence of 
these two results taken together: {\em Horn's conjecture} \cite{Horn} 
from 1962, which
gives a recursive system of inequalities, since the relevant Schubert
calculus questions can be cast as lower-dimensional
Littlewood-Richardson questions.

The main tool in this paper is the Berenstein-Zelevinsky cone \cite{BZ,Z},
and in particular the BZ polytope associated to the triple 
$(\lambda,\mu,\nu)$, in which the number of lattice points is the
corresponding Littlewood-Richardson coefficient. 
We use a new description of the BZ cone: the {\em honeycomb} model.
(The reader who is willing to grant appendix \ref{BZappendix}
does not need to absorb separately the definition of the BZ cone.)
This is a special case of a general way
of producing polyhedra that we dub {\em tinkertoy models}. This 
viewpoint gives us natural ways to interpret faces of
the BZ polytope as associated to simpler tinkertoys. In addition, the
Gel$'$fand-Cetlin system fits in this theory as associated to a
$1$-dimensional tinkertoy.

The essence of the proof is as follows. We introduce a way of indexing
(real) points in the Berenstein-Zelevinsky cone by planar pictures
called {\em honeycomb diagrams}; this identification is in appendix 
\ref{BZappendix}.
This rational polyhedral cone linearly projects to the space of triples 
$(\lambda,\mu,\nu)$ of (real) dominant weights of $\GLn$.
By \cite{BZ}, the number of integral points
(honeycomb diagrams whose vertices lie at points in the triangular lattice) 
in a fiber of this projection is the dimension
$(V_{\lambda} \tensor V_{\mu} \tensor V_{\nu})^\GLn$.
In order to work conveniently with honeycomb diagrams, 
we introduce the seemingly richer notion of a honeycomb,
since honeycombs can be seen to naturally fit into a polyhedral cone.
Then we prove the somewhat
technical theorem \ref{picturesdontlie} that honeycombs
are characterized by their diagrams 
(whose linear structure is less apparent). 

If for some large $N$ we have
$(V_{N\lambda} \tensor V_{N\mu} \tensor V_{N\nu})^\GLn > 0$, then
the fiber over $(N\lambda,N\mu,N\nu)$ of this linear projection 
contains a lattice point and is thus nonempty. By rescaling we find 
that the fiber over $(\lambda,\mu,\nu)$ is also nonempty. So the question
comes down to showing that a nonempty fiber over an integral triple
necessarily contains a lattice honeycomb. Equivalently, we want a way
of deforming a non-lattice honeycomb with integral ``boundary conditions''
$\lambda,\mu,\nu$ to a lattice honeycomb.

We do this by maximizing a linear functional,
the ``weighted perimeter'',\footnote{%
One point easily missed
is that an arbitrary choice is made in choosing this functional, making
the subsequent construction non-canonical -- but since we only seek
an existence proof, this is not a problem.}
on the polytope of honeycombs with given $\lambda,\mu,\nu$.
This picks out an extremal honeycomb,%
\footnote{%
It was not a priori obvious that the lattice point we seek occurs
as a vertex of the honeycomb polytope. In particular, not all the
vertices are at lattice points (see figure \ref{fig:noninteger}).
Nor was it plain that a single functional could be used to pick
them out uniformly for all $\lambda,\mu,\nu$. 
These facts are side consequences of the proof.}
the ``largest lift'', which we prove 
in theorem \ref{bigthm} to have very nice properties if the three weights
are suitably generic. This theorem seems to be the useful one for studying
honeycombs, and will play an equally important role in the next paper
in this series \cite{Hon2}.

It is then straightforward to prove from its nice properties
that the largest lift is integral.
A continuity argument handles the case of nongeneric triples of weights.
This ends the proof.

Not all of the framework presented in this paper is strictly necessary
if one only wishes to prove the saturation conjecture. In the very
nice paper \cite{Bu} a streamlined version of our proof is presented, 
avoiding honeycombs in favor of the hive model%
\footnote{%
The first version of this paper required the reader to absorb both
models and switched viewpoint back and forth. The paper \cite{Bu} 
was inspired by that version, and took the approach of eliminating 
honeycombs, rather than hives as is done here.}
(presented in appendix \ref{hiveAppendix}), and in particular not
requiring theorem \ref{picturesdontlie}. However, one consequence
of theorem \ref{picturesdontlie} is that honeycombs have
a very important operation called {\em overlaying}
which will be central for developments in later papers.
In the next in this series we will use the overlaying operation to study
which of Klyachko's inequalities are in fact essential \cite{Hon2}.

We thank Chris Woodward for pointing out that the saturation
conjecture gives a new proof of the weak PRV conjecture for $\GLn$,
which states that $V_{w\lambda+v\nu}$ (for $w\lambda+v\nu$ in the
positive Weyl chamber) is a constituent of $V_\lambda \tensor V_\mu$.
(This ``conjecture'' is 
nowadays known to be true for all Lie groups \cite{KM}.)
In fact one can do better, and without using saturation; 
in section \ref{sec:overlay}
there is a canonical honeycomb witnessing
each instance of the long-proven ``conjecture'', constructed by
overlaying $GL_1$-honeycombs.

We mention very briefly some connections to algebraic and symplectic geometry
(much more can be found in \cite{F,Z}). 
By Borel-Weil,
the space $(V_{\lambda} \tensor V_{\mu} \tensor V_{\nu})^\GLn$ 
is the space of invariant sections of the $(\lambda,\mu,\nu)$ line bundle
on the product of three flag manifolds.
Given two nonzero invariant sections,
one of the $(\lambda,\mu,\nu)$ line bundle and one of
the $(\lambda',\mu',\nu')$, we can tensor them together to get an
invariant section of the $(\lambda+\lambda',\mu+\mu',\nu+\nu')$ line bundle.
The geometrical fact that the flag manifold (hence the product) is
reduced and irreducible guarantees that this tensor product section is
again nonzero; this is why the set of triples $(\lambda,\mu,\nu)$ with
invariant sections forms a monoid.

This same data is involved in defining a geometric invariant theory
quotient of the product of three flag manifolds
by the diagonal action of $\GLn$.
%
%
The space $(V_{N\lambda} \tensor V_{N\mu} \tensor V_{N\nu})^\GLn$ 
is the $N$th graded piece of the coordinate ring of this quotient space.
(Klyachko's paper is a study of the semistability conditions that
arise in performing this quotient.)
The BZ counting result then says that this moduli space of triples of
flags has the same Hilbert function as a certain toric variety, and
our saturation result says that the (by definition ample) line bundle
on this moduli space actually has sections. 
W. Fulton has shown us examples in which this line bundle is not very
ample.  {\em If} one had an explicit degeneration of the moduli space
to the toric variety, one might be able to relate this
non-very-amplitude to the existence of nonintegral vertices on the
corresponding polytope of honeycombs.

The symplectic geometry connection then comes from the ``GIT quotients
are symplectic quotients'' theorem \cite{MFK}
(whose proof is essentially repeated in Klyachko's paper, 
in this special case).
In this case, 
the corresponding symplectic quotient is the space of triples of 
Hermitian matrices with spectra $\lambda$, $\mu$, and $\nu$ 
which sum to zero,
modulo the diagonal action of $U(n)$. 
$$ (\FlagsCn_\lambda \times \FlagsCn_\mu \times \FlagsCn_\nu) // \GLn
\iso \big\{ (H_\lambda,H_\mu,H_\nu) : {\rm eigen(H_\alpha) = \alpha},
H_\lambda + H_\mu + H_\nu = 0 \big\} / U(n)$$
In particular, this identification
shows directly that the existence of a $\GLn$-invariant vector in 
$V_\lambda \tensor V_\mu \tensor V_\nu$ implies the existence of a triple of
Hermitian matrices with zero sum. The reverse implication exactly
amounts to the saturation conjecture. 

B. Sturmfels has pointed out that the ``largest lift'' construction
can be interpreted as selecting a vertex of the {\em fiber polytope}
\cite{BS} of the projection from the cone of honeycombs to the
cone of triples of dominant weights. 
Combining this idea with the hypothetical degeneration of 
the moduli space to the BZ toric variety, 
this suggests that we might be able to alternately interpret
the largest lift as picking out a point in the {\em Chow} quotient 
\cite{KSZ}
of the product of three flag manifolds by the diagonal action of $\GLn$.

We thank Anders Buch, Bill Fulton, Bernd Sturmfels, Greg Warrington,
and Andrei Zelevinsky for careful readings
and many cogent suggestions; Bill we thank especially for correcting
us on a number of historical inaccuracies in the early versions.

We encourage the reader to get a feeling for honeycombs by playing 
with the honeycomb Java applet at 
\begin{center}
{\tt http://www.alumni.caltech.edu/\~{}allenk/java/honeycombs.html.}  
\end{center}

\section{Tinkertoys and the honeycomb model}\label{sec:tink}
%

We fix first a few standard notations.

A weight of $\GLn$ is a list of $n$ integers, and is dominant if the list is
weakly decreasing. So $\GLn$'s root lattice is the hyperplane of lists 
whose sum is zero. 

The ``graphs'' in this paper are rather nonstandard;
for us, a directed graph $\Gamma$ is a quadruple
$(V_\Gamma,E_\Gamma,head,tail)$ where the $head$ and $tail$ maps from
the edges $E_\Gamma$ to the vertices $V_\Gamma$ may be only {\em partially
defined} -- the edges may be semi- or even fully infinite. In particular,
any subset of the vertices and edges gives a subgraph, where the
domains of definition of the $head$ and $tail$ maps are restricted
to those edges that have their heads or tails in the subgraph.

For $B$ a real vector space, let
$Rays(B) := (B-\{0\}) / \reals_+$ denote the space of rays
coming from the origin. Each ray $d$ is in a unique line $\reals \cdot d$.
Topologically $Rays(B)$ is a sphere, $S^{\dim B -1}$.

\subsection{Tinkertoys.}
We define a {\dfn tinkertoy} $\tau$ as a triple $(B,\Gamma,d)$ 
consisting of a vector space $B$, a directed graph $\Gamma$ (possibly
with some zero- or one-ended edges), and a map $d: E_\Gamma \to Rays(B)$
assigning to each edge $e$ a ``direction" $d(e)$ in the sphere.%
\footnote{%
A related, though much more restrictive, 
definition has recently appeared in \cite{GZ},
in a context quite related to the polytope tinkertoys 
in example \ref{ex:polytope} following.
In both cases, it is sort of unnatural to fix an orientation on the
graph -- really it is the ``orientation times the direction''
that comes into play.}

\begin{Example}\label{ex:polytope}
Polytope tinkertoys.
Any polytope $P$ in $B$ gives a natural tinkertoy, just from the vertices,
the edges (oriented arbitrarily), and their directions
$d(e) := (head(e)-tail(e))/\reals^+ \in Rays(B)$.
For example, each rectangle in $\reals^2$ with edges aligned with
the coordinate axes gives us the same polytope tinkertoy
(up to isomorphism).
\end{Example}

We define a {\dfn configuration $h$ of a tinkertoy $\tau$} 
as a function $h : V_\Gamma \to B$ assigning a point of $B$ to each vertex,
such that for each two-ended edge $e$
$$ h(head(e)) - h(tail(e)) \in d(e) \cup \{\vec 0\}. $$
More generally, we say $h$ is a {\dfn virtual configuration} of $\tau$ if
for each two-ended edge $e$
$$ h(head(e)) - h(tail(e)) \in \reals \cdot d(e). $$

The set of virtual configurations is a linear subspace of the vector
space of all maps $V_\Gamma \to B$. The set of configurations is a closed
polyhedral cone in this subspace, cut out
by the conditions that that the edges be of nonnegative length;
we call it the {\dfn configuration space} or 
{\dfn cone of configurations} of the tinkertoy $\tau$.
We can use the vector space structure to define
the {\bf sum} $h_1+h_2$ of two (virtual) configurations.

If $B$ is endowed with a lattice, 
one can speak of {\dfn lattice configurations} of the tinkertoy: 
these are the ones such that
the map $h$ takes $V_\Gamma$ to lattice points in $B$.

\begin{Example}\label{ex:polyconf}
Configurations of polytope tinkertoys
(for cognoscenti of toric varieties only -- we neither use nor
prove the statements in this example).
In the case $P$ a convex lattice polytope such that the edges from
each vertex give a $\integers$-basis of the lattice,
there is an associated smooth toric variety,
and the polytope tinkertoy is just a way of
encoding the (complete) fan of the polytope.
The vector space of virtual configurations of the corresponding
polytope tinkertoy can be naturally identified with the second
equivariant cohomology group of the toric variety \cite{GZ};
the cone of actual configurations is then identified with
the equivariant K\"ahler cone, and the lattice configurations
with the equivariant Chern classes of nef line bundles.%
\footnote{%
In order to model line bundles on {\em noncompact} toric varieties 
using tinkertoys, we would need a more 
ornate definition of tinkertoy including higher-dimensional 
objects than edges, and also a more ornate definition of configuration,
assigning affine subspaces to all objects in the tinkertoy (not just
to the vertices). These additional complications only serve to obscure
the simplicity of the tinkertoys actually used in this paper.}
\end{Example}

We define a {\dfn subtinkertoy} $(B,\Delta,d|_{E_\Delta}) \leq (B,\Gamma,d)$ 
as a tinkertoy living in the same space $B$, with any subgraph 
$\Delta \leq \Gamma$, and the same assigned directions $d$ (restricted
to the subset $E_\Delta$). 
We will not have much need for morphisms of tinkertoys, but 
we do define an {\dfn isomorphism} between
two tinkertoys $(B,\Gamma_1,d_1), (B,\Gamma_2,d_2)$ in the same space $B$
as a correspondence between the two graphs, intertwining
the direction maps $d_1, d_2$.

\begin{Example}\label{ex:gc}
The Gel$'$fand-Cetlin tinkertoy.
Let $B=\reals$, $V=\{v_{i,j}\}$ for $1\leq i \leq j\leq n$, and $E$
consist of two groups of edges $\{e_{i,j}, f_{i,j}\}$, 
each $1\leq i \leq j\leq n-1$. Every edge is assigned the direction
$\reals^+$. 
$$ head(e_{i,j}) = v_{i,j+1}, \qquad tail(e_{i,j}) = v_{i,j} $$
$$ head(f_{i,j}) = v_{i,j}, \qquad tail(f_{i,j}) = v_{i+1,j+1} $$
One important subtinkertoy in this consists of the ``primary" vertices
$\{v_{i,n}\}$ and no edges.  The configurations of the
Gel$'$fand-Cetlin tinkertoy restricting to a given configuration of
the primary vertices form a polytope called the
Gel$'$fand-Cetlin polytope. Each high weight of $\GLn$ gives a (weakly
decreasing) list of integers, which we take as a lattice configuration
of the primary vertices. The lattice points in the
Gel$'$fand-Cetlin polytope are called Gel$'$fand-Cetlin patterns, and
they count the dimension of the corresponding irreducible
representation of $\GLn$. 
Note that not every configuration of the
primary vertices can be extended to a configuration of the whole
Gel$'$fand-Cetlin tinkertoy -- for this to be possible, the 
coordinates of the primary vertices {\em must} be weakly decreasing.  
\end{Example}

Return now to the general case.
If the two-ended edges $e$ of a tinkertoy $\tau$
are all of {\em positive} length in a configuration $h$, or equivalently 
$$ h(head(e))-h(tail(e)) \in d(e), $$
we call the configuration $h$ {\dfn nondegenerate}.\footnote{%
One last toric variety remark:
in the context of the configurations of polytope tinkertoys, 
nondegenerate lattice configurations correspond to ample line bundles
on the corresponding projective toric variety.}
Otherwise we say that $h$ is a {\dfn degenerate configuration}, 
each edge $e$ with $h(head(e)) = h(tail(e))$ is a {\dfn degenerate edge}
of $h$, and each vertex attached to a degenerate edge is a 
{\dfn degenerate vertex} of $h$.
Note that not every tinkertoy {\em has} a nondegenerate configuration --
for example, make a tinkertoy with vertices $x,y$ and two edges $e,f$ from
$x$ to $y$ with different given directions $d(e) \neq d(f)$. 

The following proposition is of a type standard in convex geometry:

\begin{Proposition}\label{prop:nondeg}
  If a tinkertoy $\tau$ has a nondegenerate configuration, then the
  nondegenerate configurations form the interior of the cone of
  configurations, and every configuration is a limit of nondegenerate ones.  
  If $\tau$ doesn't have any nondegenerate configurations, there is
  some edge $e$ that is degenerate in every configuration of $\tau$.
\end{Proposition}

It is worth noting that if $h$ is a degenerate configuration of a
tinkertoy $\tau$, one can associate a smaller tinkertoy $\bar\tau$
in which all the degenerate edges of $h$ have been removed, and any
two vertices connected by a series of degenerate edges have been
identified; the configuration $h$ then descends to a {\em nondegenerate}
configuration of this smaller tinkertoy. In this way each face of
the cone of configurations can be identified with the full cone of
configurations of a smaller tinkertoy.

\subsection{The relevant vector space $B$ for this paper's tinkertoys.}
From here on out, all our tinkertoys are going to live in the same
space $\B := \{(x,y,z) \in \reals^3 : x+y+z = 0\}$, a plane containing
the triangular lattice $\BZ$. This plane has three {\dfn coordinate
  directions} $(0, -1, 1), (1, 0, -1), (-1, 1, 0)$, and each direction
$d(e)$ will be one of these.

In particular, as one traverses the interval assigned to an edge
by a configuration,
one coordinate remains constant while the other two trade off,
maintaining zero sum.  We will call this the {\dfn constant coordinate} 
of the edge in the configuration.

\begin{Example}\label{ex:gl2}
The $GL_2$ honeycomb tinkertoy, in figure \ref{gl2tink}.
(This case is too small to see why these are named ``honeycombs''.)
This tinkertoy has one vertex which is the tail of three edges in
the three coordinate directions,
three vertices that are each the head of three such edges, for a total
of four vertices and nine edges (six of which have no tails
and are thus \semiinfinite).
\begin{figure}[htbp]
  \begin{center}
    \input{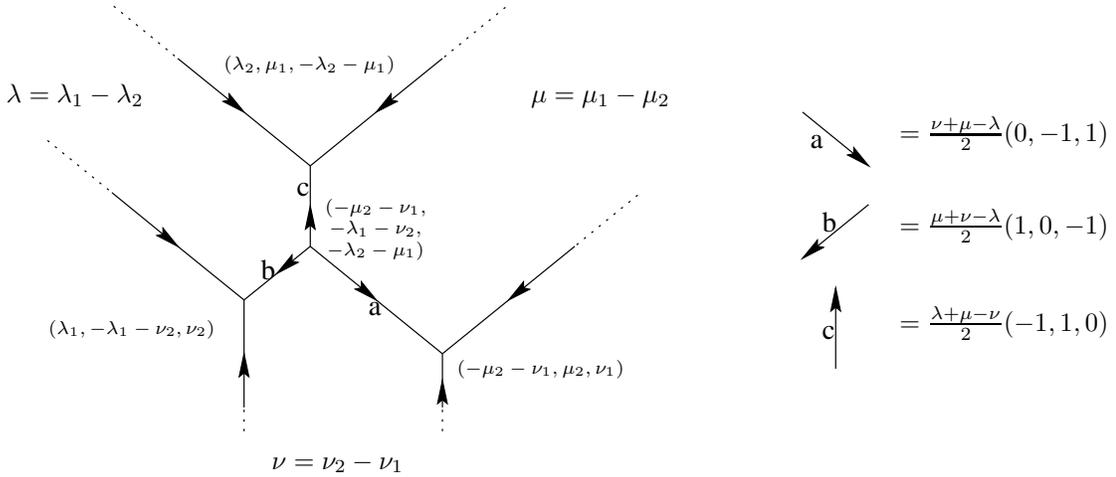}
    \caption{A configuration of the $GL_2$ honeycomb tinkertoy, 
      the vertices labeled with their coordinates in $\B$.
      The lengths of the two-ended edges are given at right,
      in terms of the separations $\lambda,\mu,\nu$.}
    \label{gl2tink}
  \end{center}
\end{figure}

The constant coordinates on the six \semiinfinite edges
determine the configuration, and
looking at the central vertex, we see that their sum must be zero.
This equality is only sufficient for the existence of a {\em virtual}
configuration: any actual configuration will satisfy also the triangle
inequalities on the separations
$\lambda = \lambda_1-\lambda_2, \mu=\mu_1-\mu_2,\nu = \nu_1-\nu_2$
between the pairs of \semiinfinite edges going going in a coordinate
direction.

There is a concise way to describe the set of integral 
coordinates $\lambda_1,\lambda_2,\mu_1,\mu_2,\nu_1,\nu_2$ 
that arise in configurations of the $GL_2$ honeycomb tinkertoy: 
they are exactly those such that the tensor product
$V_{(\lambda_1,\lambda_2)} \tensor
V_{(\mu_1,\mu_2)} \tensor
V_{(\nu_1,\nu_2)}$ of the corresponding representations of $GL_2$
contains an invariant vector. (Proof sketch: 
the requirement that the sum be zero is equivalent to asking that the
center of $GL_2$ act trivially on the tensor product. Then the triangle
inequalities are familiar from $SL_2$ theory.)
In these cases, the configuration is unique, and so too is the
invariant vector (up to scale).
\end{Example}

Much of the rest of this section is about generalizing this example
to general $GL_n$, which quite amazingly can also be performed
in the plane $\B$.


\subsection{The infinite honeycomb tinkertoy, and $GL_n$ honeycomb
tinkertoys.}

As promised, $B = \B$. Let $V$ be the set of points%
\footnote{%
This is $sl_3$'s weight lattice minus its root lattice. 
Presumably there is a deep meaning to this -- perhaps relating to the
$A_2$ web diagrams in \cite{Ku} -- but we did not uncover it.}
$$V := \{(i,j,k) \in \integers^3_{\sum=0} : 3\hbox{ doesn't divide } 2i+j\}.$$
For each vertex $(i,j,k)\in V$ such that $2i+j \equiv 2 \bmod 3$,
put on three outwardly directed edges, ending at the vertices
$(i-1,j+1,k),(i,j-1,k+1),(i+1,j,k-1)$.
These will be the vertices and edges of a directed graph $\Gamma$,
in which every edge is two-ended, and each vertex has three attached edges,
either all in or all out (depending on $2i+j \bmod 3$).

\begin{figure}[Hht]
  \begin{center}
    \leavevmode
    \epsfig{file=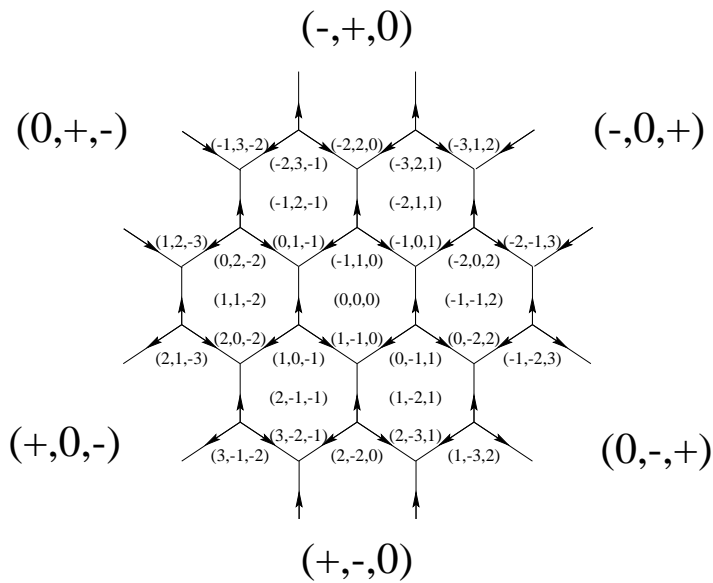,height=3in}
    \caption{A small region in the standard configuration of the 
      infinite honeycomb tinkertoy. The six adornments on the boundary
      show how the coordinates change as one moves in that direction.}
  \end{center}
\end{figure}

The {\dfn infinite honeycomb tinkertoy} is then 
$(\B,\Gamma,d)$, where the direction $d(e)$ of an edge is its direction
$(head(e)-tail(e))/\reals^+$,
and the inclusion map $V \into \B$ defines a
nondegenerate configuration 
of this tinkertoy.

\begin{figure}[htbp]
  \begin{center}
    \input{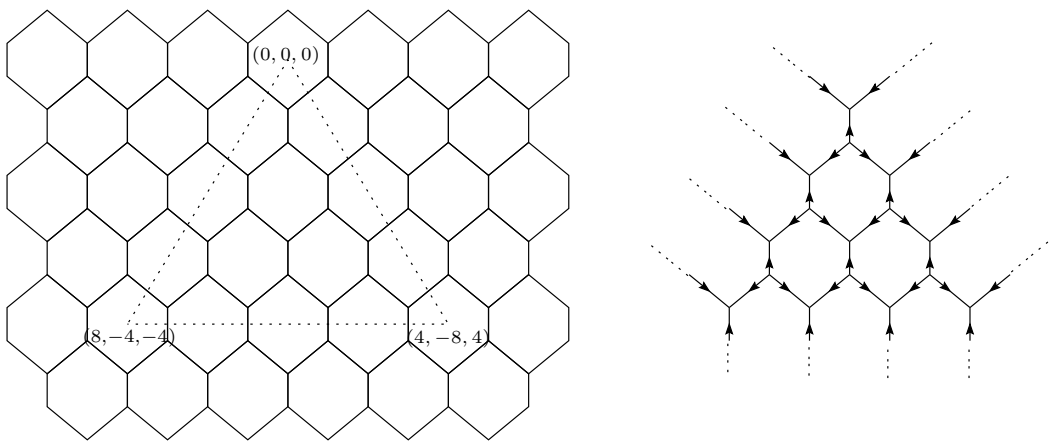}
    \caption{The triangle in the infinite honeycomb tinkertoy containing
      the $GL_4$ honeycomb tinkertoy, and that tinkertoy on its own.}
    \label{gl4tink}
  \end{center}
\end{figure}

The {\dfn $GL_n$ honeycomb tinkertoy} $\tau_n$ is the subtinkertoy of the
infinite honeycomb tinkertoy whose vertices are the $(i,j,k) \in V$
contained in the triangle $j + 3n \geq i \geq k \geq j$ (automatically
in the interior), and all their attached edges.
This tinkertoy has $3n$ tailless edges;
we will call these the {\dfn boundary} edges of this tinkertoy.
The $n=4$ example can be seen in figure \ref{gl4tink}.

There is a more general notion of honeycomb tinkertoy that we defer
until section \ref{sec:degen}.

We call a configuration $h$ of $\tau_n$ a 
{\dfn honeycomb} or {\dfn $\tau_n$-honeycomb}.
It is important to distinguish the ontological levels here --
the $GL_n$ honeycomb {\em tinkertoy} is an abstract graph with some 
labeling by directions, whereas a {\em honeycomb} 
is the additional data of an actual configuration of that tinkertoy in $\B$. 

\begin{figure}[htbp]
  \begin{center}
    \input{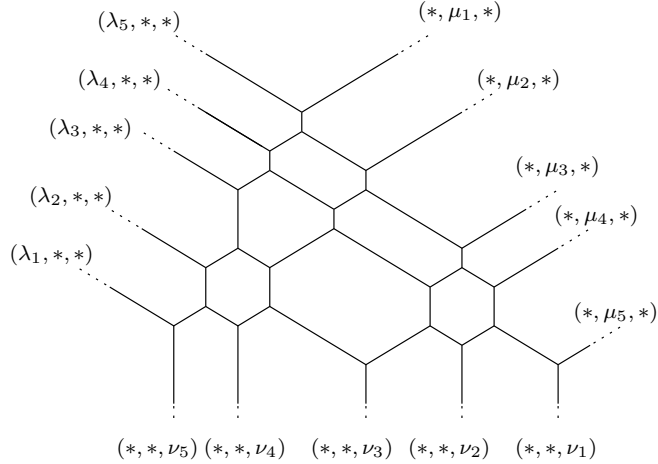}
    \caption{The constant coordinates on the boundary edges
      of a $\tau_5$-honeycomb. (The stars are the nonconstant coordinates).}
    \label{fig:LRhoneybdry}
  \end{center}
\end{figure}

Given a $\tau_n$-honeycomb $h$, we can read off the constant coordinates
on the $3n$ \semiinfinite edges starting from the southwest and
proceeding clockwise.\footnote{%
If our notion of ``configuration'' of a tinkertoy explicitly assigned
lines in $B$ to edges in $E_\Gamma$, not just to those with a head or tail, 
we could regard this as the 
restriction of a configuration of the $GL_n$ honeycomb tinkertoy to
the subtinkertoy consisting of the boundary edges and no vertices.}
Denote these $\lambda_1,\ldots,\lambda_n$, $\mu_1,\ldots,\mu_n$, 
$\nu_1,\ldots,\nu_n$, as in figure \ref{fig:LRhoneybdry}; these are
the {\dfn boundary conditions} of the honeycomb.

Let $\HONEY(\tau_n)$ denote the cone of $\tau_n$-honeycombs,
and $\BDRY(\tau_n)$ the cone of possible boundary conditions 
$(\lambda,\mu,\nu)$ of $\tau_n$-honeycombs. That is to say,
$\BDRY(\tau_n)$ is the image in $(\reals^n)^3$ of the map
$h \mapsto$ its boundary conditions $\lambda,\mu,\nu$.

Our purpose in introducing honeycombs is to calculate 
Littlewood-Richardson coefficients, the dimensions 
$\dim (V_\lambda \tensor V_\mu \tensor V_\nu)^{\GLn}$.
We do this by linearly relating $\tau_n$-honeycombs 
to Berenstein-Zelevinsky patterns%
\footnote{%
Gleizer and Postnikov have recently given
a way of relating honeycomb configurations to
Berenstein-Zelevinsky patterns that is
very different from ours \cite{GP}.}
in an appendix, where we establish the
$\integers$-linear equivalence of the Berenstein-Zelevinsky cone
with the space $\HONEY(\tau_n)$.

That equivalence has the following consequence:

\begin{Theorem*}[from appendix \ref{BZappendix}]
Let $\lambda,\mu,\nu$ be a triple of dominant weights of $\GLn$,
and $\tau_n$ the $GL_n$ honeycomb tinkertoy.
Then the number of lattice $\tau_n$-honeycombs
whose \semiinfinite edges have constant coordinates
$\lambda_1,\ldots,\lambda_n,\mu_1,\ldots,\mu_n,\nu_1,\ldots,\nu_n$
as in figure \ref{fig:LRhoneybdry}
is the Littlewood-Richardson coefficient
$\dim (V_\lambda \tensor V_\mu \tensor V_\nu)^{\GLn}$.
\end{Theorem*}

This generalizes the $GL_2$ case we did before as example \ref{ex:gl2}.

\begin{Example}\label{ex:adjsq}
In figure \ref{fig:adj} we calculate the tensor square 
of the adjoint representation of $GL_3$.
The corresponding Littlewood-Richardson rule calculation,
throwing away partitions with more than three rows, gives
$$ V_{(2,1,0)}\tensor V_{(2,1,0)} = V_{(4,2,0)} \oplus  
V_{(3,2,1)}^{\oplus 2} \oplus V_{(4,1,1)} 
\oplus V_{(3,3,0)} \oplus V_{(2,2,2)}.$$
(Recall that to turn the $S_3$-symmetric honeycomb formulation back into
a tensor product decomposition, one must reverse and negate the weight
considered the ``output''.)

\begin{figure}[Hhtb]
  \begin{center}
    \leavevmode
    \epsfig{file=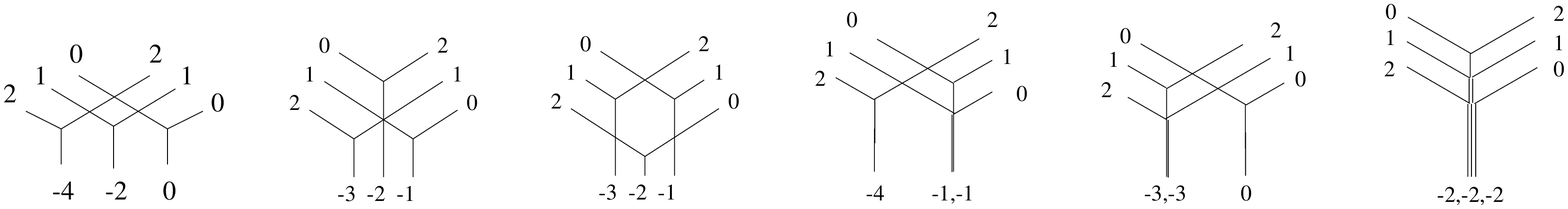,width=6in}
    \caption{The honeycombs computing the tensor square of
      $GL_3(\complexes)$'s adjoint representation. N.B. the 
      \semiinfinite edges are labeled with their constant
      coordinates -- usually we will label them with their
      multiplicities.}
    \label{fig:adj}
  \end{center}
\end{figure}
\end{Example}

One elementary consequence of this theorem is that 
the sum of the constant coordinates on the boundary edges is zero.
Proof: by linearity and continuity, 
it is enough to check on lattice $\tau_n$-honeycombs.
On the representation theory side, the sum of the constant coordinates
gives the weight of the action of the center of $\GLn$, which must
be trivial for there to be any invariant vectors. QED. We now set up
a more direct proof by a sort of Green's theorem argument, proving
some other results in tandem.

For an edge $e$ in the $GL_n$ honeycomb tinkertoy $\tau_n$, 
let the closed {\dfn interval} $I_{h,e}$ be defined by
$$
\begin{array}{ccccccc}
I_{h,e}:=
                 &\big(h(head(e)) - \reals_{\geq 0} \cdot d(e)\big) &
        \,\cap\, &\big(h(tail(e)) + \reals_{\geq 0} \cdot d(e)\big) & \\
&\bullet\!\!\longrightarrow& &\longleftarrow\!\!\bullet&
\end{array} 
$$
where if $head(e)$ or $tail(e)$ is undefined
the corresponding term is omitted.
(Note that $\tau_n$ has no zero-ended edges, so this intersection is
never over the empty set.)
Say that a curve $H$ in $\B$ {\dfn intersects the
configuration $h$ transversely} if $H$ contains no points $h(v)$, 
and intersects each interval $I_{h,e}$
transversely or not at all.

\begin{Lemma}\label{lem:jordan}
  Let $\tau_n$ be the $GL_n$ honeycomb tinkertoy,\footnote{%
This lemma applies word-for-word to the more general honeycomb
tinkertoys defined in the next section.}
 $h$ a $\tau_n$-honeycomb, 
  and $\gamma$ a piecewise-linear Jordan curve in $\B$ intersecting
  $h$ transversely. For each edge $e$ let $\gamma_e$ be the number
  of times $I_{h,e}$ pokes through $\gamma$ from the inside to
  the outside, minus the number from the outside to the inside
  (the total\footnote{%
    This is also equal to the ``number of heads of $e$ landing inside
    $\gamma$ minus the number of tails'' -- a perhaps misleading
    phrase, since each number is only ever zero or one!}
  will be $1$, $0$, or $-1$).
  
  1. The sum over $e\in E$ of the unit vectors in the
  directions $d(e)$, weighted by $\gamma_e$, is the zero vector. 

  2. The sum over $e\in E$ of the constant coordinates on $e$,
  weighted by $\gamma_e$, is zero. 
\end{Lemma}

\begin{proof}
  Since $\tau_n$ has its standard configuration,
  which is nondegenerate, by proposition \ref{prop:nondeg} the
  nondegenerate $\tau_n$-honeycombs are open dense in the cone of all
  $\tau_n$-honeycombs. The space of $\tau$-honeycombs 
  intersecting $\gamma$ transversely is open in the space
  of all $\tau$-honeycombs, and each of the functionals 
  above is obviously piecewise linear (and continuous) on it.
  Therefore the nondegenerate $\tau$-honeycombs are open dense 
  in the ones intersecting $\gamma$ transversely, and by continuity
  it suffices to check the lemma for them.
  
  If $\gamma$ encloses one (or no) vertices the statement is easily checked.
  Otherwise we can connect two points on $\gamma$ by a path within the
  interior intersecting $h$ transversely and separating the vertices
  into two smaller groups (see figure \ref{fig:jordan}).  This
  shortcut gives us two new Jordan curves, $\gamma_1$ and $\gamma_2$.
  One checks that each of the above functionals satisfies $f(\gamma) =
  f(\gamma_1) + f(\gamma_2)$.  By induction the two terms on the
  right-hand side are zero, and therefore the left is also.
  \begin{figure}[htbp]
    \begin{center}
      \epsfig{file=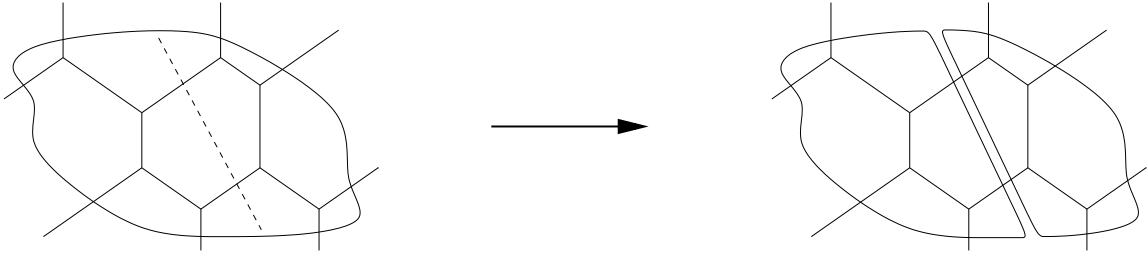,width=6in}
      \caption{Replacing a single Jordan curve by two, 
        with the dashed line as the shortcut.}
      \label{fig:jordan}
    \end{center}
  \end{figure}
\end{proof}

In particular, if we take our Jordan curve to be (a PL approximation to) 
a very big circle, we recover the previous result that the
sum of the constant coordinates on the boundary edges is zero.

\subsection{Eliding simple degeneracies.}

Recall from above that we call a tinkertoy configuration $h$ 
degenerate if some edge has length zero.
%
%
This can be regarded as a configuration of a simpler tinkertoy $\bar \tau$, 
in which the two vertices collapsed together have been identified 
(and the edge removed). The configuration map $h: V_\Gamma \to B$ descends
to give a configuration $\bar h$ of $\bar \tau$.
In this way the faces of a configuration
cone can be identified with configuration cones of simpler tinkertoys.

The case of interest to us is when a single edge of a honeycomb tinkertoy
degenerates to a point -- or more generally, when no two degenerate
edges share a vertex. 
In this very special case there is an even simpler tinkertoy to consider,
where these five edges $E,F,G,H,I$ and two vertices $x,y$ 
are not replaced by four edges $E,F,H,I$ and one vertex $x=y$, 
but two edges $E=I,F=H$ and no vertices at all. (In particular,
though $E$'s head/$I$'s tail $x$ is removed, the identified edge
$E=I$ gets the head $r$ and tail $p$; similarly $F=H$ the head $s$
and tail $q$.)
\begin{figure}[Hht]
  \begin{center}
    \leavevmode
    \epsfig{file=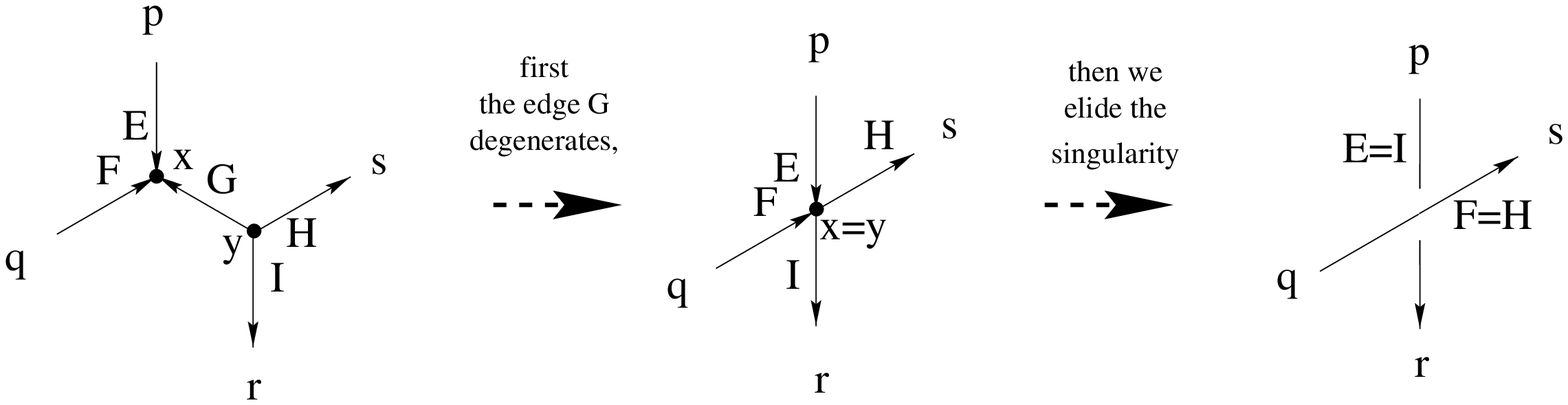,height=1.5in}
  \end{center}
\end{figure}
We will call this modification of a tinkertoy {\dfn eliding} the
edge $G$, or the vertices $x,y$, or just eliding the singularity.
A tinkertoy created by eliding a number of edges in a honeycomb tinkertoy
we call a {\dfn post-elision tinkertoy}.
Note that there is no analogue of this for a degenerate edge in a
general tinkertoy -- it is crucial that $d(E)=d(I), d(F)=d(H)$
so that the $d$ map is well-defined on the post-elision tinkertoy.
Note also that worse degenerations, 
as occur in the tensor product calculation in figure \ref{fig:adj}, 
do not usually allow the vertex to be removed. 

We will say a honeycomb $h$ has only {\dfn simple degeneracies} if no
two degenerate edges meet in a vertex.  In this case we will typically
elide the degenerate edges in the sense of the paragraph above.

It is not readily apparent what the degrees of freedom of a tinkertoy are. 
However, for post-elision tinkertoys one can say something useful.

\begin{Lemma}\label{rigid}
Let $\tau$ be the $GL_n$ honeycomb tinkertoy.\footnote{%
Again, this lemma extends word-for-word to the general honeycomb tinkertoys
defined later.}
Let $h$ be a $\tau$-honeycomb some of whose degeneracies are simple,
and $\bar\tau$ be a post-elision tinkertoy obtained by eliding some of $h$'s 
simple degeneracies, so $h$ descends to a configuration $\bar h$ of $\bar\tau$.
Let $\gamma$ be a loop (undirected) in the underlying graph of $\bar\tau$
containing only nondegenerate vertices of $\bar\tau$ (necessarily trivalent).
Then there is a one-dimensional family of configurations of $\bar\tau$,
starting from $\bar h$, in which one moves only the vertices in $\gamma$.
\end{Lemma}

\begin{proof}
%
We can assume that the loop doesn't repeat vertices; if it does,
it will have subloops that do not (in which case we will move a proper
subset of $\gamma$'s vertices).

Orient the loop, and label the vertices with signs based on whether the
loop turns left or right at the vertex, as in figure \ref{breathing}. 
Because of the angles, if all
the left-turn vertices move so as to shrink their non-loop edges by a fixed
length $\epsilon$,
whereas the right-turn vertices move so as to extend their
outgoing edge by the same $\epsilon$,
the angles remain unchanged -- i.e. we have a new configuration.
\begin{figure}[Hht]
  \begin{center}
    \leavevmode
    \epsfig{file=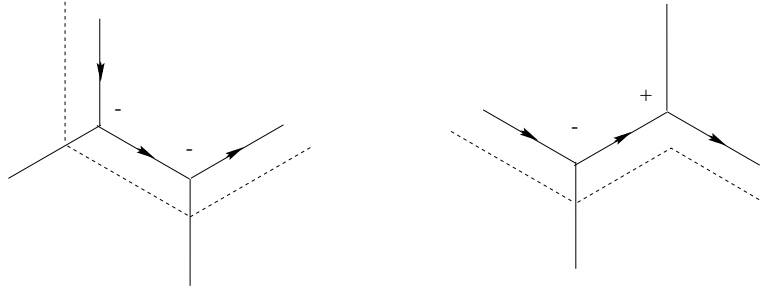,height=1.5in}
    \caption{Successive vertices turning the same direction, and in
      opposite directions, and where they could move (dashed).
      The signs indicate the change in length of the non-loop edge
      at the vertex.}
    \label{breathing}
  \end{center}
\end{figure}
\end{proof}
%

Note that no loop can go through \semiinfinite edges. So in this 
previous lemma
we're only studying degrees of freedom which leave the \semiinfinite
edges in place. Also, because we can orient the loop either way, the loop
can breathe both in and out.

%
%

\section{The diagram and degeneracy graph of a honeycomb,
and reconstructing a honeycomb from its diagram}\label{sec:degen}

\subsection{Honeycomb tinkertoys and honeycomb diagrams.}

For $h$ a configuration of a tinkertoy $\tau = (B,\Gamma,d)$, define 
the {\dfn diagram $m_h$ of the configuration $h$} to be a measure on $B$, 
the sum
$$ 
m_h := 
\sum_{e\in E_\Gamma} \hbox{Lebesgue measure on the interval $I_{h,e}$}. 
$$
This is a little more information than the set-theoretic union
$\Union_{e\in E_\Gamma} I_{h,e}$, in that it remembers 
multiplicities when edges are directly overlaid.
Note that we can recover the union from $m_h$, as its support $\supp m_h$.

In this section we will prove that a $\tau_n$-honeycomb is reconstructible
from its diagram. 
There is a stronger statement -- that every measure on $\B$
that looks enough like the diagram of a honeycomb is indeed
the diagram of a unique honeycomb (up to a trivial equivalence) --
but it requires a more general definition of honeycomb tinkertoy
than the $GL_n$ honeycomb tinkertoys we have met so far.

Define a {\dfn hexagon} in the infinite honeycomb tinkertoy
as the six vertices around a hole, i.e. 
a 6-tuple of vertices  $(i-1,j+1,k),(i,j-1,k+1),(i+1,j,k-1),
(i+1,j-1,k),(i,j+1,k-1),(i-1,j,k+1)$ where $3$ divides $2i+j$.

Define a {\dfn honeycomb tinkertoy} $\tau$ as a
subtinkertoy of the infinite honeycomb tinkertoy satisfying five
conditions:
\begin{enumerate}
\item $\tau$ is finite
\item (the underlying graph of) $\tau$ is connected
\item each vertex in $\tau$ has all three of its edges (which may
  now be one-ended)
\item $\tau$ contains a vertex (it's not just a single no-ended edge)
\item if four vertices of a hexagon are in $\tau$, all six are.
\end{enumerate}
(It is slightly unfortunate to rule out the infinite honeycomb tinkertoy
itself, but it would be more unfortunate to have to say
``{\em finite} honeycomb tinkertoy'' throughout the paper.)  
We will call the configuration of $\tau$ restricted from the defining
configuration of the infinite honeycomb tinkertoy the 
{\dfn standard configuration} of $\tau$.

\begin{figure}[htbp]
  \begin{center}
    \epsfig{file=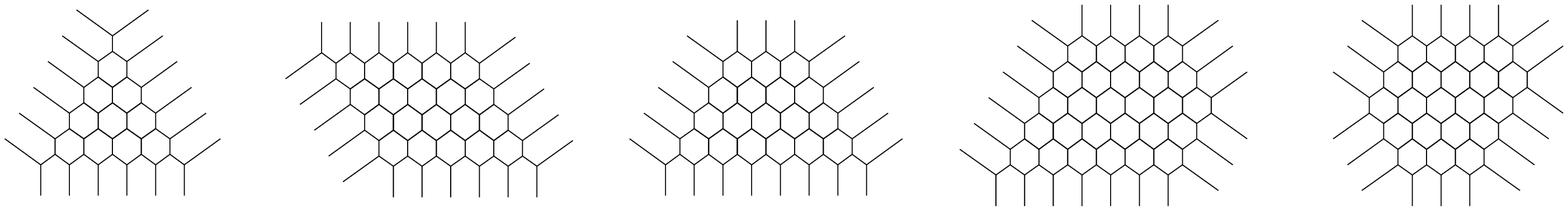,width=6in}
    \caption{The standard configurations of some honeycombs. The
      boundary edges are \semiinfinite.}
    \label{fig:standardhoneycombs}
  \end{center}
\end{figure}

A honeycomb tinkertoy has a number of \semiinfinite edges
in each of the three coordinate directions and their negatives.
Call this ordered 6-tuple, counted clockwise from North, 
the {\dfn type} of the honeycomb tinkertoy.
One fact we prove later, in lemma \ref{classifyhoneycombs},
is that any two honeycomb tinkertoys of the same type 
are isomorphic -- that the honeycomb
tinkertoys presented in figure \ref{fig:standardhoneycombs}
essentially capture all the types.
(In fact they are better than `isomorphic'; they differ only by 
translation within the infinite honeycomb tinkertoy.)
We give the more axiomatic definition above to make
it easy to check, in lemma \ref{Dregions},
that a certain subtinkertoy is itself a honeycomb tinkertoy.

If $\tau$ is a honeycomb tinkertoy, we define a {\dfn $\tau$-honeycomb} $h$
to be a configuration of $\tau$, and will speak of $h$'s type (meaning
the type of $\tau$). 
Just as in the case of the $GL_n$ honeycomb tinkertoy $\tau_n$,
we use $\HONEY(\tau)$ to denote the cone of
$\tau$-honeycombs, and $\BDRY(\tau)$ to denote the cone of possible
constant coordinates on the set of one-ended edges of $\tau$.

It is quite easy to describe the local structure of the diagram
of a (perhaps degenerate) honeycomb $h$. In the pictures to follow
of honeycomb diagrams we label edges with their multiplicities
(multiple of Lebesgue measure on the line).

\begin{Lemma}\label{lem:types}
  Let $m_h$ be the diagram of a honeycomb $h$.
  Each point $b \in \B$ has a neighborhood in which $m_h$
  satisfies one of the following:
  \begin{enumerate}
  \item $m_h = 0$
  \item $m_h$ is equal to a natural times Lebesgue measure on a
    coordinate line through $b$
  \item $m_h$ matches one of the following, 
    up to rotation (here the edge multiplicities are naturals):
  \begin{figure}[htbp]
    \begin{center}
    \epsfig{file=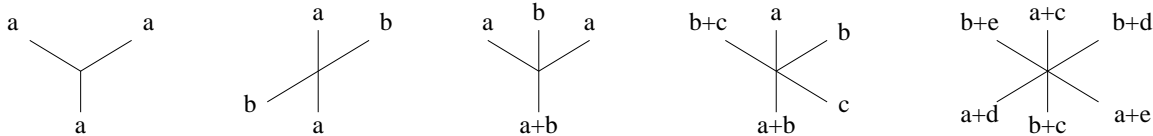,width=6in}
   \caption{The {\bf Y}, crossing, rake, $5$-valent, and $6$-valent vertices.}
    \label{fig:possvertices}
    \end{center}
  \end{figure}
  \end{enumerate}

  There are only finitely many points in the third class,
  the {\dfn vertices} of the honeycomb diagram $m_h$.
  If one thinks of each of the edges meeting a vertex as 
  pulling on the vertex with a tension equal
  to its multiplicity, these are exactly the ways for the vertex to
  experience zero total force.
\end{Lemma}

\begin{proof}
  That the pictures above are the only possibilities, and the
  finiteness of the number of points in the third class, each follow
  from the finiteness of the edge set of a honeycomb tinkertoy (so the
  neighborhood can be shrunk to avoid the $I_{h,e}$ not actually
  meeting the point $b$), and the directions $d(e)$ being multiples of
  $60^\circ$ from North.  It remains to be sure the multiplicities are
  constrained as claimed.

  By lemma \ref{lem:jordan}, the sum of the unit outgoing edges of a
  vertex (weighted by their multiplicities) must be zero. So if both
  a direction and its negative appear with positive multiplicity, we
  can subtract one from each and continue. Eventually we must get to
  a vertex of the first type, {\bf Y}, or nothing at all.
\end{proof}

Each of the points $b$ of the third class in the above we will call
a {\dfn vertex of the diagram} $m_h$. 

To see that each of these vertex types actually occurs, start with the 
honeycombs in figure     \ref{fig:standardhoneycombs}
(or larger versions with more hexagons) and degenerate all
the two-ended edges to points.
As we will show in this section, 
collapsing the honeycombs in figure \ref{fig:standardhoneycombs}
is essentially the {\em only} way to produce the singular vertices
in lemma \ref{fig:possvertices}.

\subsection{The dual graph $D(\tau)$ of a honeycomb tinkertoy $\tau$.}

Call the lattice $\{(i,j,k) \in \BZ: 3 \hbox{ divides } 2i+j\}$ 
the {\dfn root lattice} (in that it is the root lattice of $\lie{sl}_3$),
and define the {\dfn root lattice triangle} around a vertex $(i,j,k)$
of the infinite honeycomb tinkertoy to be the three points of
the root lattice at $L^1$-distance 2 from $(i,j,k)$ (these are the
three closest points). We will need the small triangular graph
made from a root lattice triangle, and also the region enclosed.

Fix a honeycomb tinkertoy $\tau$.
Define $D(\tau)$,
the {\dfn dual graph of $\tau$}, to be the union of
the root lattice triangles around the vertices in $\tau$ -- 
this has one vertex in each region in the standard configuration of $\tau$ 
(including the
unbounded ones), with an edge connecting two $D(\tau)$-vertices if the
corresponding $\tau$-regions share an edge.  (In this way edges in
$D(\tau)$ correspond to perpendicular edges in $\tau$.)  In particular
$D(\tau)$ is naturally embedded in $\B$ -- it has more structure than
just the abstract dual graph.

\begin{figure}[htbp]
  \begin{center}
    \epsfig{file=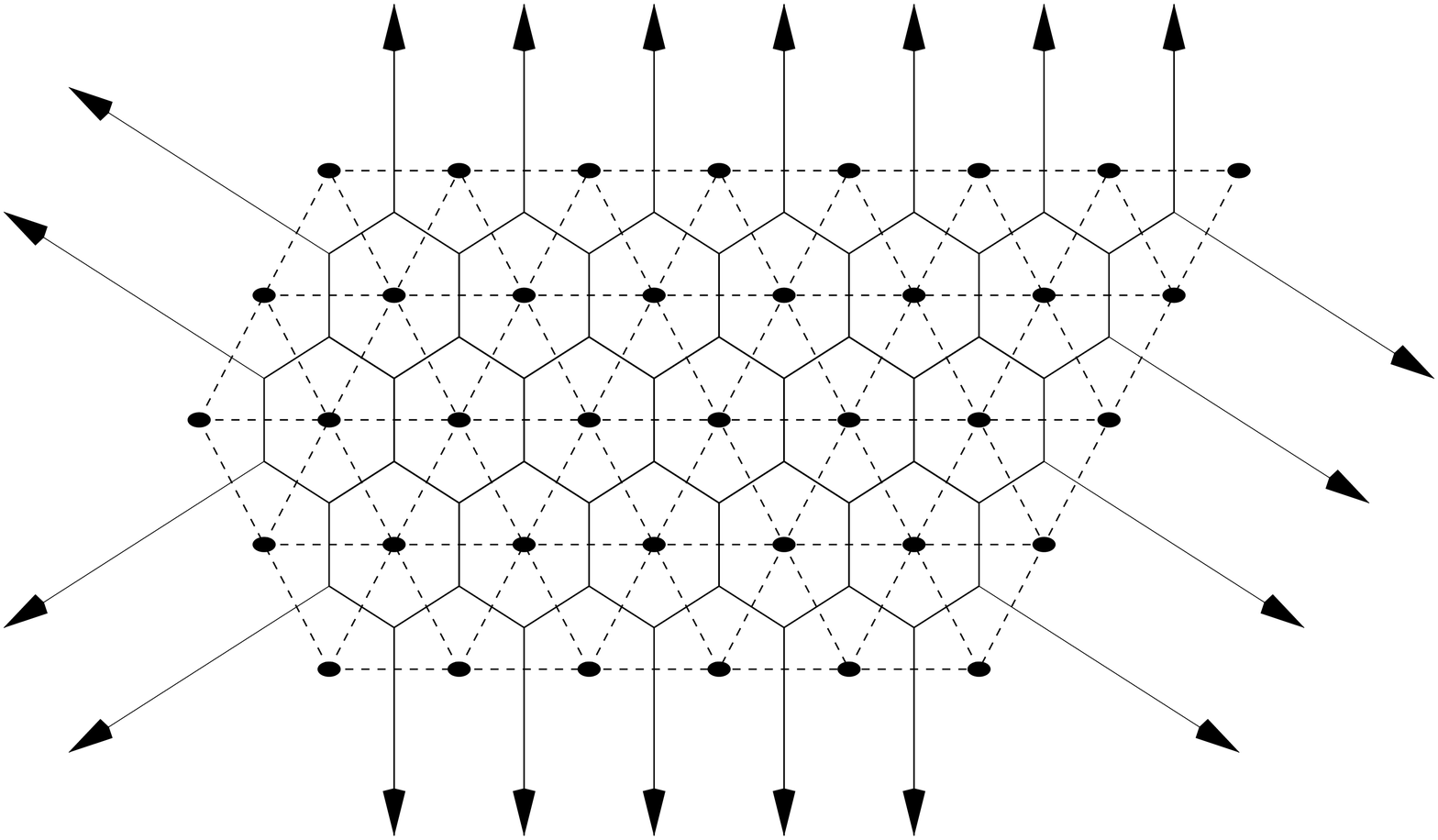,height=1.5in}
    \caption{A honeycomb tinkertoy $\tau$ of type $(7,0,4,5,2,2)$ in 
      solid lines and its corresponding $D(\tau)$ in dashed lines,
      shown superimposed in $\B$.}
    \label{fig:Dgraph}
  \end{center}
\end{figure}

\begin{Lemma}\label{Dconvex}
  Let $\tau$ be a honeycomb tinkertoy.
  Then the region bounded by its dual graph $D(\tau)$ is convex,
  and is thus characterized (up to translation) by its 6-tuple of
  edge-lengths.  Conversely, every convex union of root lattice
  triangles arises as a $D(\tau)$.
\end{Lemma}

\begin{proof}
Each vertex in $\tau$ gives us three vertices in $D(\tau)$
(by either adding 1 to, or subtracting 1 from, each of the three coordinates),
and thus a small triangle. Two connected vertices in $\tau$ share two
of these three vertices, so their corresponding triangles in $D(\tau)$
intersect in an edge (and not just a vertex).

Since $\tau$ is by assumption connected, any two of these triangles
are connected by a chain of triangles sharing common edges.
This shows that $D(\tau)$ bounds a single region, not a disconnected
set, nor two regions intersecting in only a vertex.

It remains to prove this region is convex. If not,
it has an internal angle of more than $180^\circ$
going around some boundary vertex.
This would mean four successive triangles out of the same vertex are
in $D(\tau)$. On the $\tau$ side, that means four successive vertices 
around a hexagon are in $\tau$. 
But that forces the whole hexagon to be in $\tau$.
So the vertex was not actually on the boundary, a contradiction.
This establishes the convexity of $D(\tau)$.

The converse is simple:
given a convex union $D$ of root lattice triangles we wish to realize as a
dual graph $D(\tau)$,
take $\tau$ to be the subtinkertoy of the infinite honeycomb tinkertoy
lying within $D$ (those vertices, and all their edges). 
This is easily seen to be a honeycomb tinkertoy whose
dual graph $D(\tau)$ is $D$.
\end{proof}

\begin{figure}[Hht]
  \begin{center}
    \leavevmode
    \epsfig{file=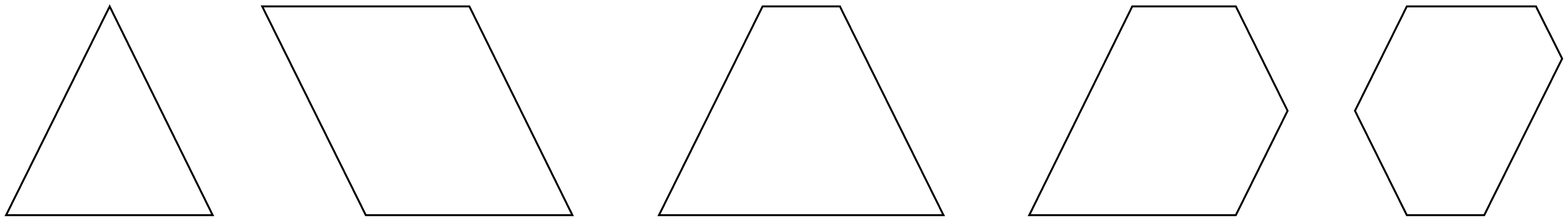,width=5in}
    \caption{The possible shapes of the region bounded
      by a dual graph $D(\tau)$.}
    \label{fig:panoply}
  \end{center}
\end{figure}

So the external angles are restricted to $0^\circ, 60^\circ,$ and 
$120^\circ$. To rotate once, the total of the external angles must be 
$360^\circ$, so there are five types, depending on the number and
ordering of the two kinds of angles; see figure \ref{fig:panoply}.

Not any $6$-tuple of edge-lengths will do; the boundary of $D(\tau)$
must be a closed curve. This gives two linear conditions that exactly
match the zero-tension property of lemma \ref{lem:types}.
So we've essentially classified the possible $D(\tau)$.

To get an equally good hold of $\tau$, we need a lemma saying we
can reconstruct $\tau$ from $D(\tau)$. That will follow from
a study of the length-minimizing paths connecting two points
in a honeycomb tinkertoy, which we dub {\dfn geodesics}.
(Note that there are a typically a great many paths with this
minimum length.)

The following lemma says that $\tau$ is ``geodesically convex''
inside the infinite honeycomb tinkertoy.

\begin{Lemma}\label{geodesics}
  Let $\tau$ be a honeycomb tinkertoy, $A,B$ two vertices of $\tau$,
  and $\gamma$ a geodesic between $A$ and $B$ 
  in the infinite honeycomb tinkertoy.
  Then the vertices of $\gamma$ are in $\tau$.
\end{Lemma}

\begin{proof}
We induct on the length of $\gamma$, assuming that the lemma is proven for
all $A$ and $B$ with geodesics shorter than $\gamma$.

Since $\tau$ is connected, $A$ and $B$ are connected under {\em some}
path $\delta$ in $\tau$.  Let $\Omega$
be the collection of hexagons enclosed by the concatenation
$\gamma + \delta$.  We can assume the cardinality of $\Omega$ is minimal 
among all possible paths $\delta$ in $\tau$ that connect $A$ and $B$. 
If $\Omega$ is empty, we are done, so suppose for contradiction that $\Omega$ 
is non-empty.

If $\gamma + \delta$ is not a Jordan curve, we can break it into smaller
pieces and use the induction and minimality hypotheses.  Hence
we may assume $\gamma + \delta$ is Jordan.

The curve $\delta$ cannot contain three consecutive edges of a hexagon
in $\Omega$; if it did, then by property 4 required 
of honeycomb tinkertoys, all the
vertices of this hexagon would be in $\tau$.  Then we could
``flip'' $\delta$ to go around the other side of the hexagon, which would
remove that hexagon from $\Omega$ and contradict the minimality assumption.
Thus we may assume that $\delta$ does not contain three consecutive
edges of any hexagon in $\Omega$.

To finish the contradiction we shall invoke

\begin{Sublemma}
  Suppose $\gamma + \delta$ is a Jordan curve whose interior $\Omega$
  is a nonempty collection of hexagons.  Suppose further that $\delta$
  does not contain three consecutive edges of a hexagon in $\Omega$.
  Then $\gamma$ is longer than $\delta$.
\end{Sublemma}
 
\begin{proof}
  Suppose for contradiction that there was a counterexample $\gamma + \delta$
  to this sublemma.  We may assume that this counterexample has a minimal
  number of hexagons in $\Omega$.  We may assume that $\Omega$ contains
  more than one hexagon, since the sublemma is clearly true otherwise.

  From hypothesis, $\delta$ does not contain three consecutive edges
  of any hexagon in $\Omega$.  We now claim that $\gamma$ also does not
  contain three consecutive edges of any hexagon in $\Omega$.  For,
  if $\gamma$ did contain three such edges, one could then ``flip'' these 
  edges across the hexagon; this would preserve the length of $\gamma$ and 
  therefore contradict the minimality of $\Omega$.  (If the flip 
  operation causes $\gamma + \delta$ to cease being Jordan, eliminate 
  redundant edges and divide into connected components).

  Now traverse $\gamma+\delta$ once in a counter-clockwise direction,
  so that $\Omega$ is always to the left.  At every vertex one turns
  60 degrees in a clockwise or counterclockwise direction.  From the
  above considerations we see that one cannot execute two consecutive 
  counter-clockwise turns while staying in the interior of $\delta$,
  since this would imply that $\delta$ contains three consecutive edges
  of a hexagon in $\Omega$.  Similarly one cannot execute two
  consecutive counter-clockwise turns while staying in the interior
  of $\gamma$.  Thus, with at most four exceptions, every 
  counter-clockwise turn in $\gamma + \delta$ is immediately followed
  by a clockwise turn.  But this contradicts the fact that we must
  turn 360 degrees counterclockwise as we traverse $\gamma + \delta$.
\end{proof}

Since $\gamma$ was assumed to be a geodesic, we have 
the desired contradiction.
\end{proof}

\begin{Lemma}\label{classifyhoneycombs}
  One can reconstruct a honeycomb tinkertoy $\tau$ 
  from its dual graph $D(\tau)$: its vertices are
  $$ \{ (i,j,k) \in \BZ : 3 \hbox{ doesn't divide } 2i+j,
  \hbox{ the root lattice triangle around $(i,j,k)$
    is in $D(\tau)$} \}.$$
  The number of \semiinfinite edges of $\tau$ in a particular
  direction is equal to the length of the corresponding edge of $D(\tau)$.
  In particular, honeycomb tinkertoys are characterized (up
  to translation in the infinite honeycomb tinkertoy) by their type.
\end{Lemma}

\begin{proof}
Each point in $\tau$ is in the set above, tautologically -- a point in $\tau$
leads to the three points in $D(\tau)$, which lead back to the same point 
being in the set above.

Now fix a triangle $T$ in $D(\tau)$; we wish to show that the center of that
triangle is necessarily in $\tau$. From a vertex in $D(\tau)$ one can
infer that at least one of the six neighboring points in 
$\{(i,j,k): i+j+k=0, 3\hbox{ doesn't divide } 2i+j\}$ is in $\tau$.
For example, the presence of the dotted (North) 
vertex in figure \ref{triangleT}
says that up to left-right reflection, there must be a $\tau$ vertex in 
one of the regions labeled $1$, $2$, $3$ or $4$.

\begin{figure}[Hht]
  \begin{center}
    \leavevmode
    \epsfig{file=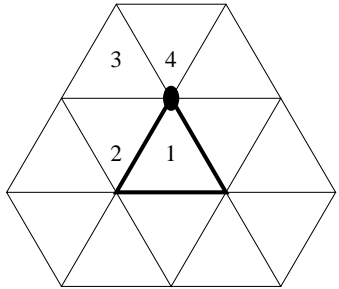,height=1.5in}
    \caption{Possibilities for the $\tau$ vertex causing 
      the North $D(\tau)$ vertex in this triangle.}
    \label{triangleT}
  \end{center}
\end{figure}

If there's a $\tau$ vertex in region $1$, we're done. If there's a $\tau$
vertex in region $2$ or $3$, and another $\tau$ vertex producing the 
existence of the Southeast vertex of $D(\tau)$, we can find (in figure
\ref{triangleT234}) a geodesic in the infinite honeycomb connecting the two
that goes through the center. Then by lemma \ref{geodesics} about
such geodesics, the center is necessarily a vertex of $\tau$.

\begin{figure}[Hht]
  \begin{center}
    \leavevmode
    \epsfig{file=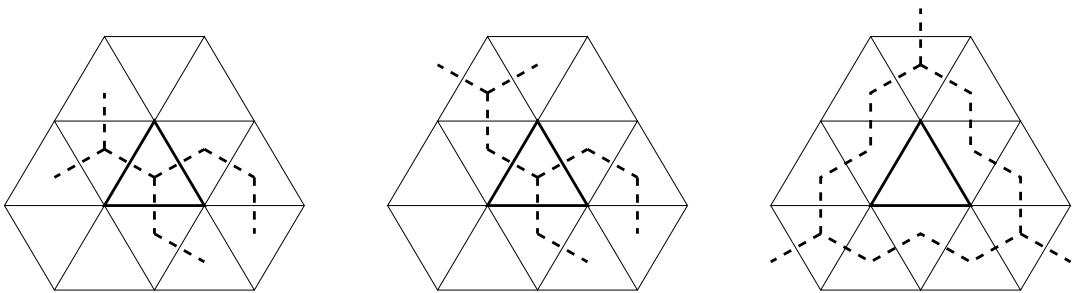,height=1.5in}
    \caption{Geodesics through the center in cases $2,3$, 
      and through $2,3$ in case $4$.}
    \label{triangleT234}
  \end{center}
\end{figure}

The remaining case occurs when the only $\tau$ vertex can be found in
region $4$, in all three rotations of the diagram. But then by connecting
two of those $\tau$ vertices with geodesics we find $\tau$ vertices in
regions $2$ and $3$, reducing to the previous case.
\end{proof}

\subsection{The degeneracy graph $D(h)$ of a honeycomb $h$.}

Fix a honeycomb tinkertoy $\tau$, and $D(\tau)$ its dual graph.
For $h$ a $\tau$-honeycomb,
let $D(h)$ be the subgraph of $D(\tau)$ with the same 
vertex set but an edge between two vertices 
only if the corresponding (perpendicular) edge of $h$ 
is nonzero. This we will call the {\dfn degeneracy graph} 
of the honeycomb $h$ (and is also embedded in $\B$).

\begin{figure}[Hht]
  \begin{center}
    \leavevmode
    \epsfig{file=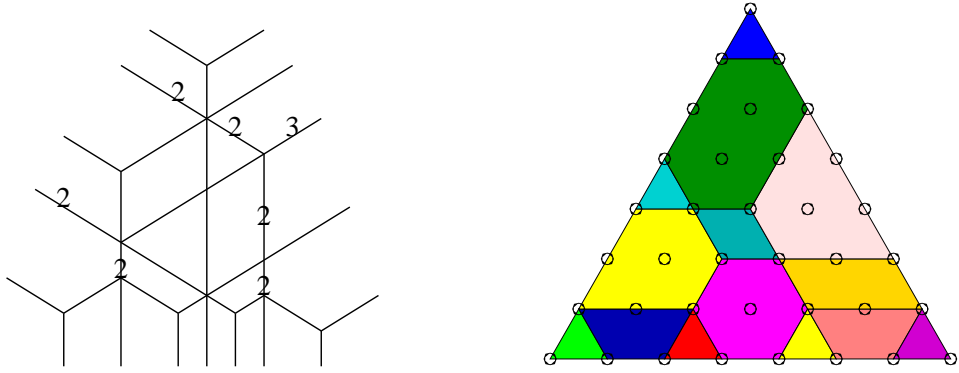,width=5in}
    \caption{A very degenerate honeycomb $h$, and the corresponding $D(h)$.}
    \label{fig:Dhregions}
  \end{center}
\end{figure}

Our goal is to classify the regions in $D(h)$, and establish that they
correspond to vertices in the diagram of $h$, thereby classifying the
possible vertices in the diagram and their preimages under $h$
(\'a la figure \ref{fig:standardhoneycombs}).  

\begin{Lemma}\label{Dregions}
  The regions in the degeneracy graph $D(h)$ of a honeycomb $h$ are convex.
\end{Lemma}

\begin{proof}
Fix a region $\Omega$ in $D(h)$, and choose a vertex $x$ of $\tau$
which lies in the closure of $\Omega$.  Then its image $h(x)$
is in the support of the diagram of $h$.

Consider the subtinkertoy of $\tau$ consisting of those vertices in $\tau$
which map to $h(x)$ under $h$, together with their associated edges.
Let $\sigma$ be the connected component%
\footnote{Actually, we shall see in the next lemma 
that the tinkertoy necessarily has only one connected component.}
of this tinkertoy that contains $x$.
Then $\sigma$ is a honeycomb tinkertoy in its own right: the only
non-trivial observation required is that if four vertices of a hexagon
map to $h(x)$, then all six vertices must map to $h(x)$.

Chasing down the definitions we see that the region bounded by $D(\sigma)$
is just $\Omega$.  The claim then follows from lemma \ref{Dconvex}.
\end{proof}

\begin{Lemma}\label{hDhduality}
  The regions in the degeneracy graph $D(h)$ of a honeycomb $h$ 
   correspond to the vertices in the diagram $m_h$ of $h$.
\end{Lemma}

\begin{proof}
One direction is clear: if two vertices in $\tau$ give 
root lattice triangles in 
the same region of $D(h)$, those vertices are connected by a series of
degenerate edges in $h$, 
and therefore have collapsed to the same vertex of the diagram of $h$.
So for each region $X$ in $D(h)$, we can speak of $X$'s $h$-vertex.
The converse is more difficult; we must show that two distinct regions
$X$ and $Y$ in $D(h)$ give vertices in the diagram of $h$ that
are physically separated. 

One case is easy. If $X$ and $Y$ are adjoining regions, then the 
existence of the edge in $D(h)$
separating them says that the corresponding edge in $h$ is nonzero.
But this is exactly the displacement between $X$'s $h$-vertex and
$Y$'s $h$-vertex, so they are in different places, as desired.
For nonadjoining $X$ and $Y$ we will have to add together a
bunch of such nonzero displacements and hope to get a nonzero sum.

Let $x$ be a generic point in $X$, and $y$ a generic point in $Y$,
and let $\overline{xy}$ the straight-line path connecting them.
Then $\overline{xy}$ does not intersect the vertices of $D(h)$
and only intersects the edges of $D(h)$ transversely.

\begin{figure}[Hht]
  \begin{center}
    \leavevmode
    \epsfig{file=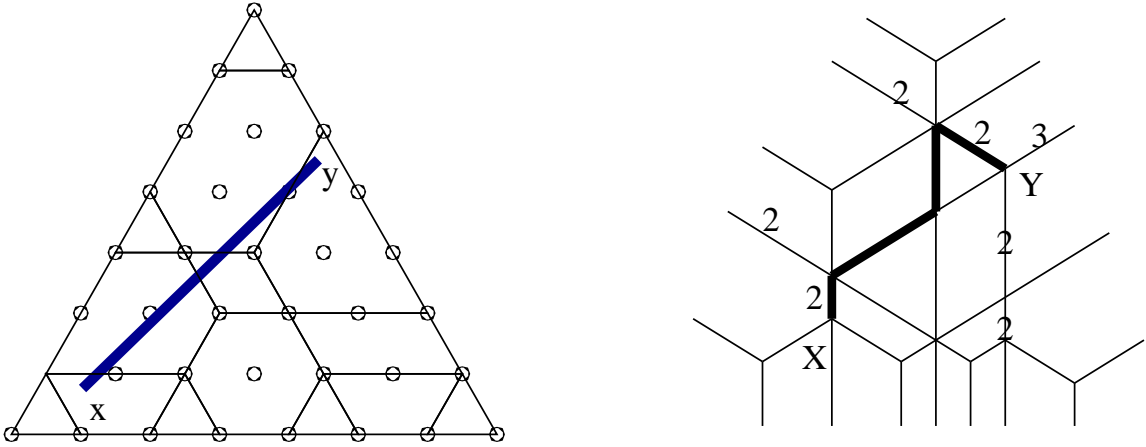,width=6in}
    \caption{A straight line from region $X$ to region $Y$ in $D(h)$, 
      and the corresponding path in $h$.}
    \label{fromXtoY}
  \end{center}
\end{figure}

We can now compute the vector in $\B$ separating the points in $h$
corresponding to the regions $X$ and $Y$.
Each time $\overline{xy}$ crosses a wall -- an edge still left in $D(h)$ --
there is an associated nonzero displacement, the vector difference
of the $h$-vertices of the regions on the two sides of the wall.
Adding up all these displacements we get the total vector difference we seek,
the displacement of the vertices corresponding to $X$ and $Y$.
(It is worth emphasizing that this
gives a path in $h$ itself, as indicated in figure \ref{fromXtoY}.)
But note now that each individual
term has positive dot product with the vector $y-x$, 
since it is perpendicular to the wall that $\overline{xy}$ has just
crossed through (and in the correct direction).

So therefore the whole sum has positive dot product with $y-x$, 
so is nonzero, and therefore $X$'s $h$-vertex and $Y$'s $h$-vertex
are not in the same place in the diagram.
\end{proof}

In particular, since the preimage of a vertex of a honeycomb is
connected, it is a honeycomb tinkertoy in its own right.

There is a succinct way to sum up the results of this section.
Define an {\dfn abstract honeycomb diagram}\footnote{%
This concept matches the {\em web functions} of \cite{GP};
we prefer to avoid this terminology, though, for fear of 
confusion with the intriguingly similar `webs' of \cite{Ku}.
A. Postnikov of \cite{GP} has informed us that he also knew 
theorem \ref{picturesdontlie}.}
as a measure $m$ on $\B$ such that 
\begin{enumerate}
\item in a neighborhood of each point $b \in \B$, 
  $m$ is a nonnegative {\em real} linear combination of the Lebesgue measures
  on the six coordinate rays out of $b$, satisfying the zero-tension
  property of lemma \ref{lem:types}
\item only finitely many points, which we naturally call the 
  {\dfn vertices} of $m$, use more than two rays in this combination
\item $\supp m$ is not a collection of parallel lines.
\end{enumerate}

By lemma \ref{lem:types}, the diagram of a honeycomb is an abstract
honeycomb diagram. 
The following theorem states that, up to a trivial ambiguity,
every abstract honeycomb diagram with integral edge-multiplicities
is the diagram of a unique honeycomb.

\begin{Theorem}\label{picturesdontlie}
  Let $m$ be an
  abstract honeycomb diagram with integral edge-multiplicities.
  Then $m$ is the diagram of a configuration $h$ of
  a honeycomb tinkertoy $\tau$, where $\tau$ is uniquely determined by
  $m$ up to unique isomorphism; and given $\tau$ the honeycomb $h$ is unique.
\end{Theorem}

\begin{proof}
  We break the diagram-to-honeycomb reconstruction into two steps:
  from the diagram to the degeneracy graph, 
  then from the degeneracy graph to the honeycomb.

  {\em Finding the degeneracy graph $D(h)$.}
  We first need to determine the dual graph $D(\tau)$ containing $D(h)$.
  Since each of the vertices in $m$ satisfies the zero-tension
  property of lemma \ref{lem:types}, the whole diagram does,
  by the same Green's theorem argument as in lemma \ref{lem:jordan}.
  Correspondingly, there does exist a convex lattice region
  whose edge-lengths are the numbers of \semiinfinite edges,
  unique up to translation. 

  By lemma \ref{hDhduality}, the vertices in $m$ are supposed to 
  correspond to the regions in $D(h)$, and we can determine the shapes
  of those regions from lemma \ref{classifyhoneycombs}. 
  It remains to fit them all together into $D(\tau)$.

  For each vertex $p$ in $m$, pick a path in
  $\supp m$ whose last vertex connects to a \semiinfinite edge. 
  Each vertex along the path corresponds to a region in $D(h)$ whose shape we
  can determine from the vertex and lemma \ref{classifyhoneycombs}.
  These regions
  glue together along edges perpendicular to the steps in the path.
  We can determine where the region
  corresponding to the last vertex in the path sits in $D(h)$:
  we know it's on a boundary edge of $D(h)$ (the boundary corresponding 
  to the direction of the \semiinfinite edge), and we know where it
  sits on that edge, by counting how many \semiinfinite edges (with
  multiplicity) going in that direction are to the right and left
  of this bunch.  Having done that, by gluing the other regions to it
  we have determined where they all sit, including that of the
  original vertex.

  The only worry then is that different paths to the boundary may suggest
  different places to locate $p$'s region inside $D(h)$. One checks that
  there is no monodromy in going around a small loop in the embedded
  graph $\supp m$ --
  this is because the hexagon-closure condition of a small path is the same
  as the zero-tension condition of lemma \ref{lem:types} --
  and therefore none in any loop. 

  {\em Finding the honeycomb tinkertoy $\tau$ and the honeycomb $h$.}
  From $D(\tau)$, using lemma \ref{classifyhoneycombs} we can construct 
  $\tau$ uniquely. Had $D(\tau)$ been chosen differently, $\tau$ would
  only change by translation within the infinite honeycomb tinkertoy.
  
  We've determined $D(h)$'s breakup into regions, and which vertex in
  $p$ corresponds to which region in $D(h)$.  But this determines the
  honeycomb $h$ -- for each vertex $v$ of $p$, take the vertices of
  $\tau$ in the interior of the $D(h)$-region corresponding to $v$ and
  map them to $v$ under $h$.
\end{proof}

The degeneracy graph $D(h)$ is thus a way of recording only 
what we might call the combinatorial information about a honeycomb $h$,
not the actual positions.
(In fact there is a tighter connection: see appendix \ref{hiveAppendix}.)
The vertices of one correspond to regions in the other, and the length
of an edge of $D(h)$ is equal to the multiplicity of the corresponding
(perpendicular) edge in $h$.\footnote{%
If the edges of $D(h)$ are given formal multiplicities
equal to the length of the corresponding edge in $h$,
the graph $D(h)$ becomes a honeycomb diagram itself 
(ignoring here the problem of dealing with the \semiinfinite edges).
This is essentially the duality in \cite{GP} on BZ triangles; see also
the remarks at the end of the next section.}

This theorem allows us to work pretty interchangeably with honeycombs
vs. honeycomb diagrams. In particular we can strengthen the result of
lemma \ref{rigid}, which only gave us a family of configurations of
post-elision tinkertoys $\bar\tau$ made from honeycomb tinkertoys $\tau$, 
to actually give us configurations of $\tau$ itself.

\begin{Corollary}
Let $h$ be a $\tau$-honeycomb,
and $\bar \tau$ be a tinkertoy obtained by eliding some of $h$'s 
degenerate edges, so $h$ descends to a configuration $\bar h$ of $\bar\tau$.
Let $\gamma$ be a loop (undirected) in the underlying graph of $\bar\tau$
passing only through nondegenerate vertices of $\bar h$.
Then there is a one-dimensional family of configurations of $\tau$,
starting from $h$, in which one moves only the vertices in $\gamma$.
\end{Corollary}

\begin{proof}
  Apply lemma \ref{rigid} to get a family of configurations of $\bar\tau$,
  and apply theorem \ref{picturesdontlie} to their diagrams;
  this produces configurations of $\tau$.
\end{proof}

\begin{figure}[Hht]
  \begin{center}
    \leavevmode
    \epsfig{file=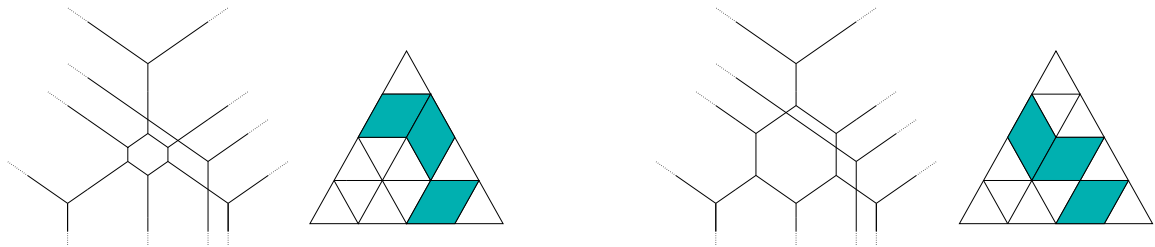,width=6in}
    \caption{Two members of a bent family of $GL_4$ honeycombs and their
      degeneracy graphs, resulting from breathing a hexagon in their
      common post-elision tinkertoy. Since their degeneracy graphs are
      different, the honeycombs lie on different faces of the cone of
      honeycombs.}
    \label{breathebend}
  \end{center}
\end{figure}

This is most interesting when breathing the loop in and out causes a
vertex of the post-elision tinkertoy $\bar\tau$ to move across an edge, 
as in figure \ref{breathebend},
something not possible for the honest honeycomb tinkertoy $\tau$.
In this case, the application of theorem \ref{picturesdontlie}
to the $1$-dimensional family of $\bar\tau$ configuration diagrams 
produces a {\em piecewise-linear} $1$-dimensional family of 
$\tau$-honeycombs that bends around the cone $\HONEY(\tau)$.

\section{Overlaying honeycombs, and the PRV conjecture}\label{sec:overlay}

The reconstruction theorem of the last section 
lets us a define a remarkable operation on honeycombs. 
(This section is not used elsewhere in this first paper.)

\begin{Corollary}[to theorem \ref{picturesdontlie}]
  Let $h$ and $h'$ be two honeycombs (perhaps of different types).
  Then (up to translation in the infinite honeycomb tinkertoy) 
  there exists uniquely a honeycomb whose diagram is the
  sum of the diagrams of $h$ and $h'$ (as measures; one adds 
  multiplicities when the edges of $h$ and $h'$ lie fully on
  top of one another).
\end{Corollary}

\begin{proof}
  Let $m$ be the sum as a measure of the diagrams of $h$ and $h'$.
  One checks straightforwardly that $m$ is an abstract honeycomb diagram
  with integral edge-multiplicities, except for connectedness.

  If $\supp m$ is not connected,
  there is a connected component $C$ of $\B \setminus \supp m$
  bounded by two different components of $\supp m$.
  The region $C$ is not convex and bounded, or else its boundary would
  have just one (polygonal) component. If $C$ is not convex
  it must be nonconvex at one of its vertices, violating 
  lemma \ref{lem:types}'s zero-tension
  condition on $m$'s vertices, contradiction. 
  
  So $C$ is convex and unbounded, hence its limit points in the circle
  at infinity $Rays(\B)$ must either be an interval or two opposite points. 
  If this is an interval, $C$'s boundary is connected, contradiction.
  So $C$ contains a line, and its boundary in $\supp m$ must be
  two parallel lines. By the zero-tension condition on honeycomb vertices,
  there can be no vertices on these lines.
  Removing them, and repeating the argument, we find that $h$ and $h'$
  are each unions of parallel lines, contrary to assumption.

  So we have the assumed connectedness, and $m$ is an abstract
  honeycomb diagram. Then by theorem \ref{picturesdontlie} it is
  the diagram of an essentially unique honeycomb.
\end{proof}

The naturality of this {\em piecewise-linear} operation on honeycombs is, 
to our minds, one of the principal advantages over the BZ formulations,
and will be central in the next paper \cite{Hon2}. (In this paper the
overlay notion is only used in a sort of local way -- the elision
operation on simple degeneracies.)

{\em Application: the weak PRV conjecture for $\GLn$.}
The so-called weak PRV conjecture (now proven in general \cite{KM}) 
states that 
if $w\lambda + v\mu$ is in the positive Weyl chamber for some Weyl
group elements $w,v$, then $V_{w\lambda + v\mu}$ is a constituent
of the tensor product $V_\lambda \tensor V_\mu$. We prove this 
for $\GLn$ as follows. 
For $i \in \{1,\ldots,n\}$,
let $h_i$ be a $GL_1$-honeycomb, i.e. a single vertex 
with three \semiinfinite edges coming off, whose coordinates are 
$(\lambda_{w(i)}, \mu_{v(i)}, -\lambda_{w(i)}-\mu_{v(i)})$.
Then overlaying all the $\{ h_i \}$, we get a lattice $GL_n$-honeycomb
with boundary conditions $\lambda$, $\mu$, and $-(w\lambda + v\mu)$.
This lattice honeycomb is a witness to 
this instance of the weak PRV conjecture.

\begin{figure}[Hhtb]
  \begin{center}
    \leavevmode
    \input{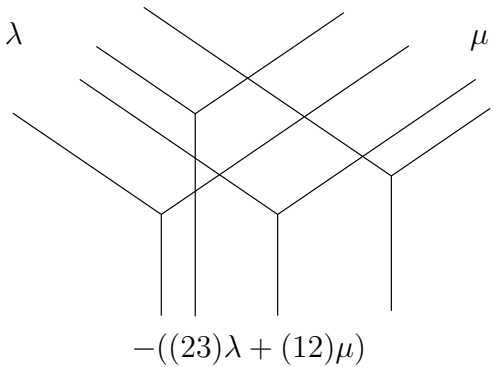}
    \caption{An example of a witness to the PRV conjecture.}
    \label{fig:PRVex}
  \end{center}
\end{figure}

The cost of working in the honeycomb model is that the linear
structure -- the fact that one can add two BZ patterns of the same
size and get another -- is geometrically a slightly more mysterious
operation on honeycombs. As we will see in a later paper in this series, 
these two operations ``add'' and ``overlay'' are intertwined by the
duality operation on honeycombs from [GP], which 
should be seen as a tropical version of the Fourier transform
relating ``times'' and ``convolve''.
On honeycombs this operation essentially amounts to replacing a
honeycomb $h$ with its degeneracy graph $D(h)$, where each edge of
$D(h)$ is given a formal ``multiplicity'' equal to the length of the
corresponding edge in $h$ (thus completing the duality of the two
graphs).  We leave the additional details of how to handle the
\semiinfinite edges to the later paper.

\section{The largest-lift map $\BDRY(\tau)\to\HONEY(\tau)$}\label{uppershell}

\begin{figure}[htbp]
  \begin{center}
    \epsfig{file=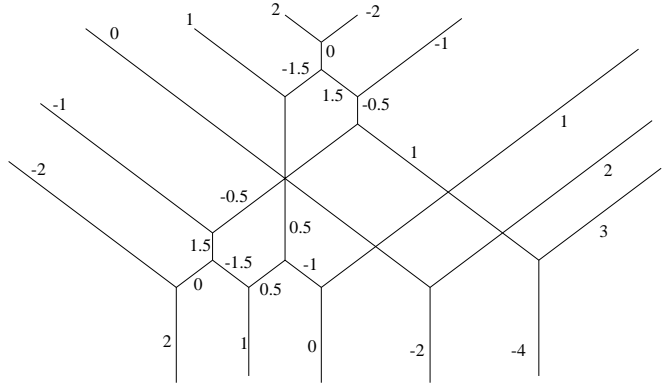,height=2in}
    \caption{A nonintegral vertex of a honeycomb polytope, with edges (all 
      multiplicity one) labeled with their constant coordinates.}
    \label{fig:noninteger}
  \end{center}
\end{figure}

In figure \ref{fig:noninteger} we see a honeycomb with integral
boundary that is not itself integral. It is a little more difficult to
see that it is in fact an extremal point on the polytope of honeycombs
with this boundary, i.e., the constant coordinates on the interior edges 
are uniquely determined unless one un-degenerates some zero-length
edges. The reader can check this by using the fact that the diagram is
in $\B$ to determine the constant coordinates on all the edges except
those in the figure eight; then the fact that the line passes through
the node on the figure eight ties down the rest of the coordinates.

It will turn out that, in a sense, 
the $6$-valent vertex is to blame for this honeycomb's bad behavior.  
In this section we develop the machinery to find extremal honeycombs
with better vertices (and thus better behavior) than this one.

Let $P$ be a (possibly unbounded) polyhedron in a vector space $U$,
$\pi: U \onto V$ a projection that restricts to a proper map $P \onto Q$, 
and $\vec w$ a generic functional on $U$.  Define the 
{\dfn ``largest lift" map} $l:Q \to P$, taking $q$ to the point in 
$P \cap \pi^{-1}(q)$ with greatest pairing with $\vec w$.  By properness,
there is a point maximizing this pairing; by genericity of $\vec w$,
this point is unique.  So the map is well-defined, and in fact is
continuous and piecewise-linear (\cite{Zi}).
In the case of $P$ a tortoise, 
$Q$ its shadow on the ground at noon,
and $\vec w$ measuring the height off the ground, 
the largest lift of a point in the shadow is the corresponding point
on the tortoise's upper shell.

Fix a honeycomb tinkertoy $\tau$.
We are interested in this largest-lift map in the case of the projection 
$\HONEY(\tau) \to \BDRY(\tau)$, which forgets the location of the vertices and
finite edges of a honeycomb, remembering only the constant
coordinates on the boundary edges.

\begin{Proposition}\label{proper}
  Let $\tau$ be the $GL_n$ honeycomb tinkertoy.
  Then the map $\HONEY(\tau) \to \BDRY(\tau)$ is proper. 
\end{Proposition}

\begin{proof}
This is guaranteed by the correspondence with the BZ cone
(in appendix \ref{BZappendix}).
\end{proof}

In fact this map is improper only if $\tau$ has 
\semiinfinite edges in all six directions, but we will not need
this fact in this paper.

The functional $wperim : \HONEY(\tau) \to \reals$
is chosen to be a generically weighted sum of the perimeters of the
(possibly degenerate) hexagons in the honeycomb, the weighting having
a certain ``superharmonicity'' property.
More exactly, let $w$ assign a real number to each of the regions in
the tinkertoy $\tau$ (vertices of $D(\tau)$), with the properties that

\begin{enumerate}
\item For each unbounded region $r$ on the exterior, $w(r)=0$
\item for each hexagon $\alpha$ surrounded by regions $\alpha_i$,
$w(\alpha) > \frac{1}{6} \sum_i w(\alpha_i)$
\item $w$ is chosen generic subject to these constraints.
\end{enumerate}

(One nongeneric such $w$ can be defined on any $D(\tau)$ by
$w(i,j,k) = -i^2-j^2-k^2$. But the set of $w$ is open, so we
can perturb this one slightly to get a generic $w$.)

Then we define the {\dfn weighted perimeter} of a $\tau$-honeycomb $h$ as
$$ wperim(h) = \sum_{\hbox{hexagons $\alpha$}} w(\alpha) 
        \hbox{ perimeter}(\alpha). $$
Since $wperim$ is defined in terms of the perimeter, it's a linear
functional on $\HONEY(\tau)$.

\begin{Lemma}\label{doinflate}
  Let $h$ be a honeycomb in which some hexagon can inflate (moving the
  vertices of the loop, as in lemma \ref{rigid}).  Then inflating
  it increases $wperim(h)$.
\end{Lemma}

  \begin{figure}[htbp]
    \begin{center}
      \epsfig{file=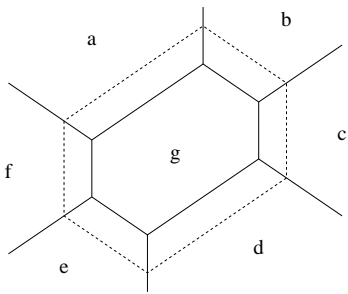,height=1.5in}
      \caption{An inflating hexagon, with each region $\alpha$
        labeled with its weight $w(\alpha)$.}
      \label{fig:inflate}
    \end{center}
  \end{figure}

\begin{proof}
  Inflating the hexagon by distance $\epsilon$ increases its perimeter
  by $6\epsilon$, while decreasing that of each of its neighbors by
  $\epsilon$. The change in $wperim(h)$ is
$$ 6g\epsilon - (a+b+c+d+e+f)\epsilon = (6g - (a+b+c+d+e+f))\epsilon > 0. $$
\end{proof}

This lemma is best understood in terms of $\HONEY(\tau)$'s linear structure,
which while transparent in the original BZ formulations, or the
hive model in appendix \ref{hiveAppendix},
is unfortunately rather opaque in the honeycomb model.
For each region $R$ in a honeycomb tinkertoy $\tau$, define the
{\dfn inflation virtual configuration} 
$\vec i_R$ which places the vertices of $R$
around the origin in $\B$ (by translating the hexagon containing $R$'s
vertices to the hexagon in the infinite honeycomb tinkertoy 
around the origin), 
but all other vertices {\em at} the origin.
Then the statement ``the hexagon $H$ can inflate in $h$, but gets stuck at
a distance $s$'' is equivalent to ``$h+\epsilon \vec i_H$ is in the cone
$\HONEY(\tau)$ for $\epsilon \in [0,s]$, but past $s$ this ray leaves
the cone''. 
And lemma \ref{doinflate} 
above is exactly the statement (if we extend $wperim$
linearly to the vector space of virtual configurations of $\tau$) that
$wperim(\vec i_R)>0$.

In a particularly degenerate honeycomb it may be impossible to
inflate any one hexagon; only certain combinations may be possible.
This next slightly technical lemma shows that certain local
changes to a honeycomb, which {\dfn molt}\footnote{%
From Webster's:
molt : to shed hair, feathers, shell, 
   horns, or an outer layer periodically  : to cast off (an outer covering) 
   periodically; specif : to throw  off (the old cuticle)}
the degeneracy, can be obtained by inflating several hexagons simultaneously. 

  \begin{figure}[htbp]
    \begin{center}
      \epsfig{file=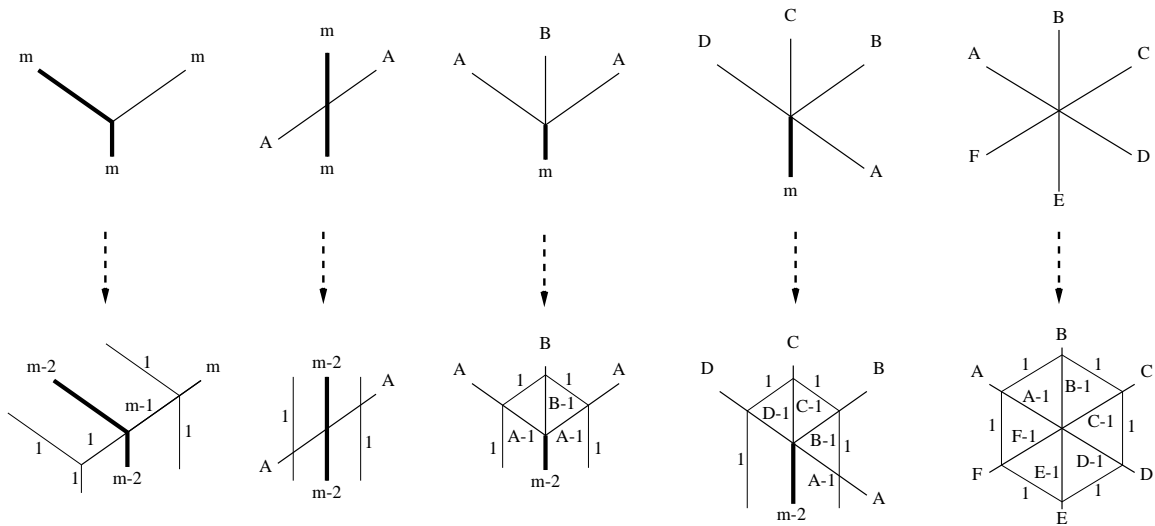,width=6in}
      \caption{How to ``molt'' a degenerate vertex.}
      \label{fig:moltingall}
    \end{center}
  \end{figure}

\begin{Lemma}\label{lem:molt}
  The virtual configurations associated to the
  recipes in figure \ref{fig:moltingall} for molting a degenerate
  vertex are in each case the sum of a set of inflation virtual
  configurations $\vec i_R$ associated to inflating a certain
  collection of regions $R$.  More precisely, each molting recipe
  inflates equally the completely degenerate hexagons, plus the
  $4$-sided regions on the sides corresponding to thick edges.
\end{Lemma}

\begin{proof}
  Since these degenerate vertices involve, by definition, a number of
  hexagons that have collapsed (to lines and even to points), it is
  rather difficult to see which hexagons must be simultaneously
  inflated to make the vertex molt. We will use the linear structure to
  get around this as follows; first add $\epsilon$ times the standard
  configuration of the honeycomb tinkertoy $\tau$. Now there is no
  degeneracy and we can point out which regions to inflate.
  Add the corresponding virtual configurations $\vec i_R$ to this
  configuration. Then subtract $\epsilon$ times
  the standard configuration and see that we do get the molted configuration.

  Note in particular that if we mark two adjacent regions for the same
  amount of inflation, their common edge doesn't move at all. 
  So it is simple to see whether an edge moves under simultaneous
  inflation, since each edge is on the boundary of exactly two
  regions; it moves if exactly one of them inflates, away from that one.

  In each of the following pictures we show 
  \begin{enumerate}
  \item the vertex (edges labeled with their multiplicities)
  \item the standard configuration of (a small example of) the
    underlying tinkertoy, with certain regions labeled in gray
  \item the result of inflating those regions some distance
  \item the same result, with the standard configuration subtracted off.
  \end{enumerate}

{\em Molting a {\bf Y} vertex} (figure \ref{fig:moltYhexes}).
We determined in lemma \ref{classifyhoneycombs} 
that {\bf Y} vertices result from the
collapse of a $(0,n,0,n,0,n)$-type honeycomb tinkertoy. 
To perform the sort of molt we want here, 
we inflate all the regions {\em except} one corner and the opposite side.

\begin{figure}[htbp]
  \begin{center}
    \epsfig{file=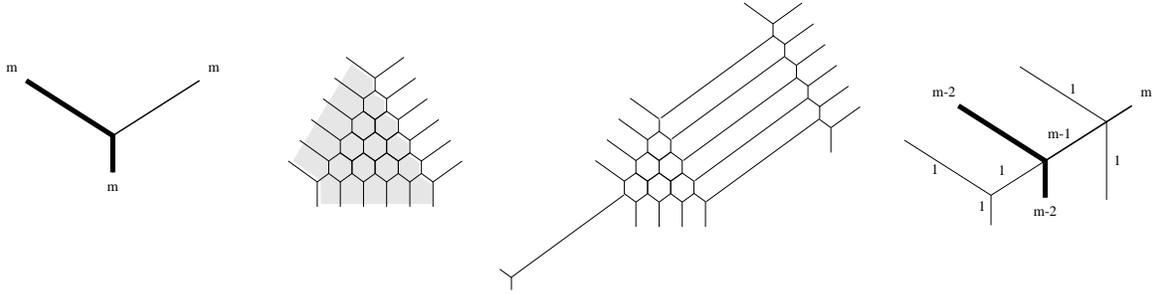,width=6in}
    \caption{From left to right: a {\bf Y} vertex; that vertex expanded
      with some regions marked for inflation (in gray); the result
      of inflation; then de-expanded again.}
    \label{fig:moltYhexes}
  \end{center}
\end{figure}

{\em Molting a crossing vertex} (figure \ref{fig:moltXhexes}).
In this case we inflate all regions except those on the left and right side.
Again, we are using lemma \ref{classifyhoneycombs} to know precisely
what is hiding in the degenerate vertex, as we will again in the remaining
cases.

\begin{figure}[htbp]
  \begin{center}
    \epsfig{file=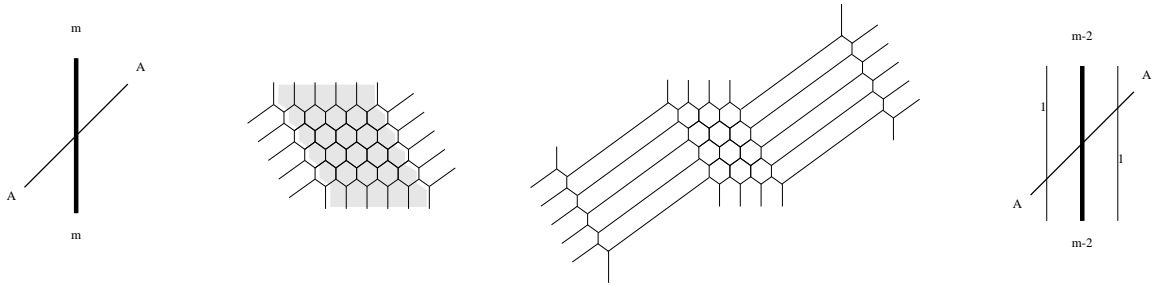,width=6in}
    \caption{Molting a crossing vertex by inflating certain regions.}
    \label{fig:moltXhexes}
  \end{center}
\end{figure}

{\em Molting a rake vertex} (figure \ref{fig:moltrakehexes}).
In this case we inflate all regions except those on the left, right, and top.

\begin{figure}[htbp]
  \begin{center}
    \epsfig{file=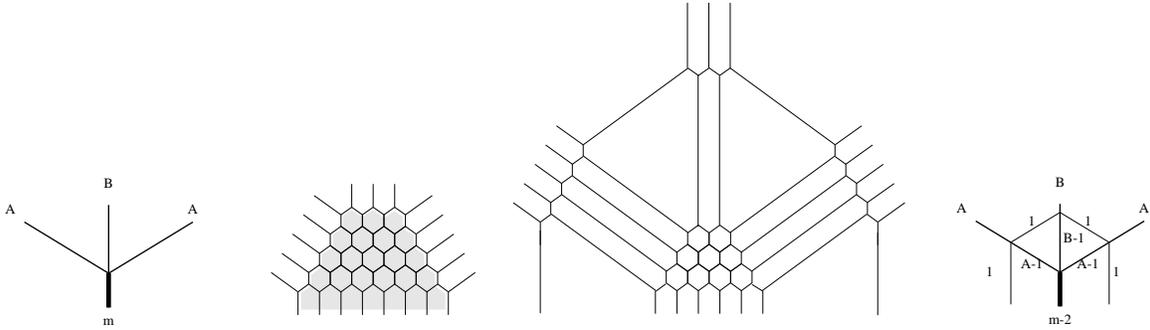,width=6in}
    \caption{Molting a rake vertex by inflating certain regions.}
    \label{fig:moltrakehexes}
  \end{center}
\end{figure}

{\em Molting a $5$-valent or $6$-valent vertex.} 
This is more of the same, so we do not take space for the pictures
(see figure \ref{fig:standardhoneycombs} for the standard configurations
of the tinkertoys).
In both the $5$- and $6$-valent case, we mark all the hexagons for 
inflation, and in the $5$-valent case we also mark for inflation the $4$-sided
unbounded regions on the side with $m$ external edges.
\end{proof}

Note that we make no statement here about the change in the 
weighted perimeter, since (except in the $6$-valent vertex case) 
some of the regions we are inflating are
unbounded, and lemma \ref{doinflate} does not pertain. 
In the theorem to follow we will only be inflating hexagons.

For $\tau$ a honeycomb tinkertoy, call $\beta \in \BDRY(\tau)$ a set of
{\dfn regular} boundary conditions if no two \semiinfinite edges of
$\tau$ in the same direction are assigned the same constant coordinate.
This terminology is taken from the $GL_n$-honeycomb case, where such
boundary conditions correspond to triples of regular dominant weights
(i.e. in the interior of the positive Weyl chamber).

The main result in this paper is the following.

\begin{Theorem}\label{bigthm}
Let $\tau$ be a honeycomb tinkertoy 
such that the map $\HONEY(\tau) \to \BDRY(\tau)$ is proper,
as in proposition \ref{proper}.
Let $\beta$ be a regular 
point in $\BDRY(\tau)$, and $l(\beta)$ the largest-lift honeycomb
lying over it (relative to a generic choice $w$).
Then $l(\beta)$ has only simple degeneracies, and if one elides them,
the graph underlying the resulting tinkertoy is acyclic.
\end{Theorem}

\begin{proof}
Since the map is assumed proper,  the concept of a ``largest lift''
makes sense, and we can go on to study its properties.
In order, we will show 
\begin{enumerate}
\item a largest lift never has $6$-valent vertices
\item a largest lift of a regular point has no edges of multiplicity $>1$
  (and therefore no rakes or $5$-valent vertices)
\item a simply degenerate largest lift has no cycles.
\end{enumerate}

Item $\# 1$ is immediate from lemma \ref{lem:molt}: if $l(\beta)$ has a
$6$-valent vertex, we can molt it by inflating a certain set of hexagons.
But by lemma \ref{doinflate} that increases the weighted perimeter.
So $l(\beta)$ was not a largest lift, contrary to assumption.

\begin{Example}
In the calculation of the tensor square of $GL_3(\complexes)$'s
adjoint representation in figure \ref{fig:adj}, we found that
two copies of the adjoint representation (tensor the determinant) appear. 
One of them has a
$6$-valent vertex and is thus not a largest lift. When that vertex
``molts'' as explained above (and the resulting hexagon is 
maximally inflated), one obtains the other honeycomb with the
same boundary, which is a largest lift.
\end{Example}

%
Proof of $\# 2$. 
Let $m$ be the maximum edge-multiplicity that appears in $h$; 
assume $m>1$ or else we're done.
Let $\Gamma$ be the subgraph of $\supp m_h$ of the edges
with multiplicity $m$ (and their vertices). By the assumption that
$\beta$ is regular, this contains none of the \semiinfinite
edges; it is bounded.
Let $x$ be a vertex on the boundary of the convex hull of $\Gamma$;
$x$ is necessarily a rake or a $5$-valent vertex.

Build a path $\gamma$ in $\Gamma$ starting at $x$, with first edge $e$,
as follows.
Declare the direction $60^\circ$ clockwise of $e$ to be the ``forbidden''
direction, and $30^\circ$ counterclockwise to be ``windward''.

Now traverse edges, coming to new vertices, going through crossings,
turning at {\bf Y}s (but not into the forbidden direction), stopping 
when you reach another rake or $5$-valent vertex.
(Conceivably one might continue through a $5$-valent vertex, if it is
lucky enough to have two edges of multiplicity $m$, but we don't do this.)
Once you start building this path from $x,e$, there are no choices,
and each step carries us a positive distance in the windward direction.
So the path doesn't self-intersect, and since the graph is bounded, 
this algorithm must terminate.
By the assumption that $m$ was the maximum edge-multiplicity, we only
come into a rake or $5$-valent along the edge labeled $m$ in
lemma \ref{lem:molt} (up to rotation and reflection).

We now attempt to simultaneously inflate all the hexagons that have
completely degenerated to the vertices along $\gamma$, plus those that
have collapsed to the edges connecting two vertices. Comparing this to
the recipes in lemma \ref{lem:molt}, we see that this exactly molts
all the vertices, and is thus a legal combination (adding a small
multiple doesn't carry us out of the cone $\HONEY(\tau)$).  
This operation inflates some hexagons, and therefore by lemma
\ref{doinflate} increases the weighted perimeter, violating the
largest-lift assumption as before.

We give an example of this in figure \ref{globalmolt}, where an entire
path $\gamma$ molts. This illustrates how the recipes for molting at a
vertex exactly fit together to give a well-defined operation on the
honeycomb. 

\begin{figure}[htbp]
  \begin{center}
    \epsfig{file=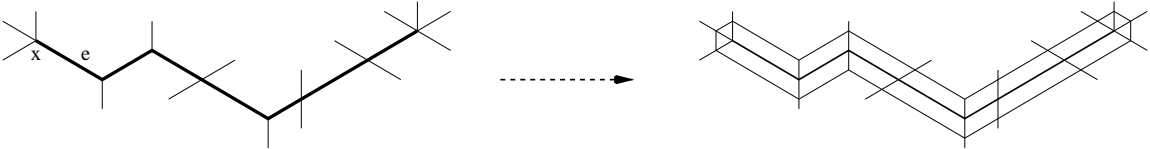,width=6in}
    \caption{Replacing a multiplicity-$m$ path with a multiplicity-$(m$-$2)$ 
path while molting a multiplicity-$1$ skin.}
    \label{globalmolt}
  \end{center}
\end{figure}

Proof of $\# 3$. 
Now that we know that our largest lift of a regular triple
only has simple degeneracies, we can apply lemma \ref{rigid}
(or rather, its strengthened version in the corollary to 
theorem \ref{picturesdontlie})
to say that any cycle in the graph underlying the post-elision tinkertoy 
gives a degree of freedom.
But by the assumed genericity of $w$, our honeycomb should be at a 
{\em vertex} of the polytope of honeycombs lying over $\beta$,
and thus have no degrees of freedom.
Hence there are no cycles in the post-elision tinkertoy.
\end{proof}

Honeycombs with nonsimple degeneracies can be seen 
in figure \ref{fig:adj} -- but only when
the bottom edge has some edges lying on top of one another,
or the honeycomb is not a largest lift.

\section{Proof of the saturation conjecture}\label{sec:sat}

We prove a general honeycomb version of the saturation conjecture.
Then we derive the actual representation theory 
saturation conjecture from its truth for $GL_n$-honeycomb tinkertoys.

\begin{Theorem}\label{saturation}
  Let $\tau$ be a honeycomb tinkertoy such that the projection
  $\HONEY(\tau) \to \BDRY(\tau)$ is proper, 
  and $w$ a generic weighting function on the regions
  satisfying the properties required in section \ref{uppershell}.
  Then the largest-lift map $\BDRY(\tau)\to \HONEY(\tau)$ is a
  piecewise $\integers$-linear map.  Consequently, any point
  in $\BDRY(\tau)$ assigning integer constant coordinates to the
  boundary edges can be extended to a lattice honeycomb.
\end{Theorem}

\begin{proof}
We already know that the largest-lift map 
is continuous, and linear on chambers. We will show by studying 
regular points in $\BDRY(\tau)$
that each of these  linear maps has integer 
coefficients. Then for any point in $\BDRY(\tau)$, even nonregular, 
we can pick a chamber of which that point is on the boundary, and show that 
over that point there lies a lattice honeycomb.

If $\beta$ is a regular configuration in $\BDRY(\tau)$,
then by theorem \ref{bigthm} from section \ref{uppershell}, 
the largest lift honeycomb $l(\beta)$ has only simple degeneracies
(its vertices only look like {\bf Y}s or crossing vertices, 
with edge-multiplicity $1$ everywhere).

By the ``elision'' construction in section \ref{sec:tink}, we can
regard this as a nondegenerate configuration of a simpler 
post-elision tinkertoy,
where each crossing point is removed, and the five edges (one
of length zero) replaced by the two lines going through.

In this tinkertoy, each vertex is degree $3$, touching some finite
and some \semiinfinite edges.
Consider the subgraph of finite edges, also acyclic, and inductively
pull off vertices of degree $1$. 
Each such vertex is connected to two \semiinfinite edges, 
whose constant coordinates determine the location of the vertex.
In particular the constant term on the finite edge coming 
out is integrally determined by those on the two \semiinfinite edges --
it is minus their sum.

\begin{figure}[Hht]
  \begin{center}
    \leavevmode
    \epsfig{file=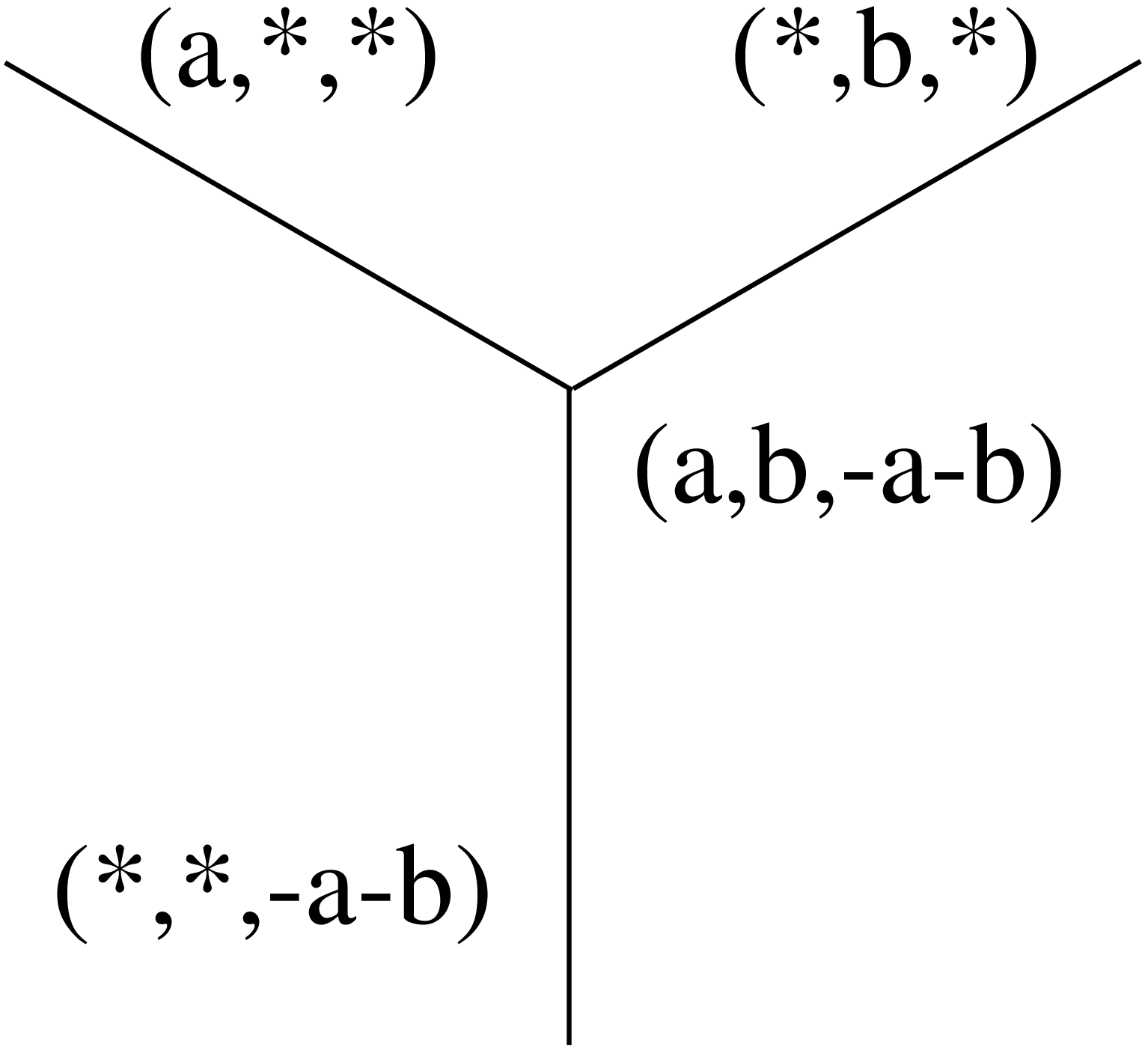,height=1.5in}
  \end{center}
\end{figure}

So we can remove the two \semiinfinite edges and the vertex, 
promoting the remaining edge to \semiinfinite, 
and recurse. Eventually all the coordinates are integrally 
determined from those on the original \semiinfinite edges.
\end{proof}

(We invite the reader to see how this argument fails on the
honeycomb in figure \ref{fig:noninteger}.)

It may be worth noting that
one can make these integral formulae very explicit. To determine the
constant coordinate on a two-ended edge $e$ somewhere
in the middle of this honey-forest $l(h)$,
let $e_{\bf Y}$ be the set of \semiinfinite 
edges $f$ such there is a path
(necessarily unique) from $f$ to $e$ going through $e$'s 
{\em right-side-up} {\bf Y} vertex. 
(Unless the post-elision tinkertoy is disconnected, the unique path from $f$
to $e$ will go either through $e$'s {\bf Y} vertex or $e$'s upside-down 
{\bf Y} vertex.) Then the constant coordinate on $e$ is the sum
of the constant coordinates on the outgoing boundary edges in $e_{\bf Y}$,
minus the corresponding sum on the incoming edges.

\begin{Corollary}[the saturation conjecture]
Let $(\lambda,\mu,\nu)$ be a triple of dominant integral 
weights of $\GLn$ such that for some $N>0$, the tensor product
$V_{N\lambda} \tensor V_{N\mu} \tensor V_{N\nu}$
has a $\GLn$-invariant vector. 
Then already
$V_{\lambda} \tensor V_{\mu} \tensor V_{\nu}$
has a $\GLn$-invariant vector.
\end{Corollary}

\begin{proof}
  Let $\tau_n$ be the $GL_n$ honeycomb tinkertoy.
  By theorem \ref{BZequiv} from appendix \ref{BZappendix},
  the fiber of the boundary-conditions map 
  $\HONEY(\tau_n) \to (\reals^n)^3$ over the point
  $(N\lambda,N\mu,N\nu)$ contains a lattice point, and is
  therefore nonempty. Therefore the fiber over $(\lambda,\mu,\nu)$
  is also nonempty, since it just the original fiber rescaled by $1/N$.
  By proposition \ref{proper} we can apply theorem \ref{saturation}, 
  which says that this fiber contains a lattice point.
  Using theorem \ref{BZequiv} again, we find that 
  $V_{\lambda} \tensor V_{\mu} \tensor V_{\nu}$
  has a $\GLn$-invariant vector.
\end{proof}

There is an analogous saturation conjecture for the tensor product of any 
number of representations, which for our purposes we can state as follows.
Let $\{\lambda_i\}_{i=1\ldots m}$ 
be a collection of dominant weights such that for some large $N$, 
the tensor product $\Tensor_i V_{N \lambda_i}$ has a $\GLn$-invariant vector. 
Then the same is true when $N$ is replaced by $1$. 

\begin{figure}[Hhtb]
  \begin{center}
    \leavevmode
    \input{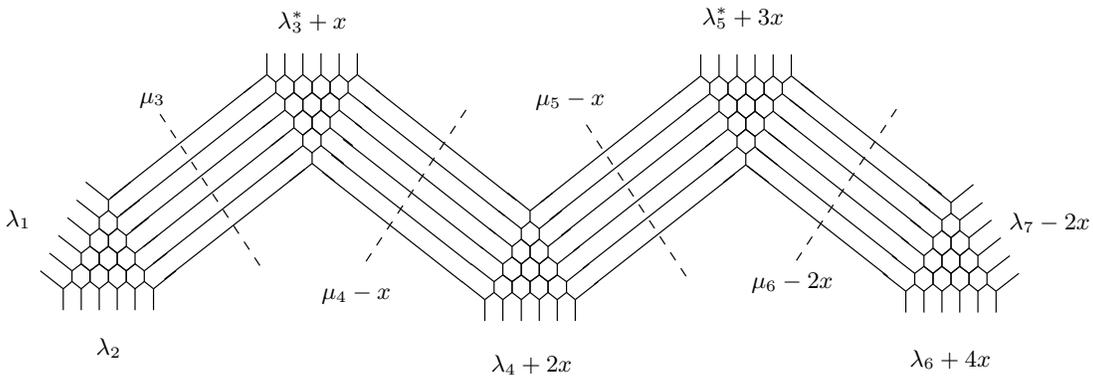}
    \caption{A honeycomb tinkertoy whose largest lift shows where to find an
      invariant vector in a $7$-fold tensor product.}
    \label{fig:satchm7}
  \end{center}
\end{figure}

An earlier version of this paper had a technically unpleasant proof
of the general saturation conjecture.
We omit the details, because since writing this paper, 
Andrei Zelevinsky has shown us how to
derive the general saturation result from the case already proven,
via standard arguments with the Littlewood-Richardson rule \cite{Z}.
The basic idea of the omitted proof is indicated in 
figure \ref{fig:satchm7}, which shows a honeycomb tinkertoy 
whose largest lifts give
witnesses to saturation in the case of a seven-fold tensor product.

As long as the inputs $\lambda_i$ are pulled apart 
sufficiently ($x\to\infty$), 
a configuration of this big honeycomb tinkertoy corresponds 1:1 to a 
`coherent' set of $(0,n,0,n,0,n)$ honeycombs (they are
far enough apart to necessarily not overlap when glued together 
into a configuration of the big honeycomb tinkertoy). 
That and repeated application of the isomorphism
$$(A\tensor B)^\GLn \iso \sum_\mu (A\tensor V_\mu)^\GLn \tensor
(V_\mu^* \tensor B)^\GLn $$
(where $\mu$ runs over all irreducible
representations of $\GLn$) let us locate a lattice point in the
analogous Berenstein-Zelevinsky polytope.

\section{A saturation conjecture for other groups}

One can phrase a na\"\i ve saturation conjecture for other groups $G$:
for any triple $(\lambda,\mu,\nu)$ of dominant weights for $G$,
if there exists a number $N$ such that
$$ (V_{N\lambda} \tensor V_{N\mu} \tensor V_{N\nu})^G > 0,
\qquad \hbox{then} \qquad
(V_{\lambda} \tensor V_{\mu} \tensor V_{\nu})^G > 0.$$
However, this conjecture is false, with counterexamples reported in
\cite{E}. A clue is provided by the fact that it is already 
false for $G=\SLn$!

Of course, the $\SLn$ case is not so different from the $\GLn$ case,
where saturation holds, and so is easily patched up; we must ask also that
the representation
$\lambda+\mu+\nu$ of the torus annihilate the center of $\SLn$.
(This is no longer a linear condition on $\lambda+\mu+\nu$, 
as it was for $\GLn$, 
exactly because $\SLn$'s center is not connected.) 

\begin{Conjecture}
  Let $G$ be a connected complex semisimple Lie group with maximal torus $T$,
  $(\lambda,\mu,\nu)$ a triple of dominant weights for $G$,
  and $N$ a positive number such that 
  $$ (V_{N\lambda} \tensor V_{N\mu} \tensor V_{N\nu})^G > 0.$$
  Then if $\lambda+\mu+\nu$ annihilates all elements of $T$ with 
  semisimple centralizer,
  $$ (V_{\lambda} \tensor V_{\mu} \tensor V_{\nu})^G > 0.$$
\end{Conjecture}

We now explain the geometric motivation for this conjecture.
Since our only goal is to make the conjecture plausible
we do not waste space on full proofs of the statements made here.

The space
$ (V_{\lambda} \tensor V_{\mu} \tensor V_{\nu})^G$ can be thought of
as the space of sections of a sheaf on the GIT quotient $(G/B)^3 // G$
(as explained in the introduction, in the case of $G = \GLn$).
However, this sheaf may not necessarily be a line bundle; it can have
singularities that get worse at orbifold points of the quotient,
and any section is required to vanish at the singularities.
The condition that $\lambda+\mu+\nu$ annihilate the center exactly
guarantees that this sheaf be a line bundle generically, so is 
certainly necessary for the existence of nonvanishing sections.

However, if the sheaf has singularities along which any section
must vanish, it stands to reason that global sections are less likely
to exist. The condition in the conjecture is exactly equivalent
to asking that the sheaf be a line bundle globally, therefore to have 
a better chance to have sections.

To be sure, this conjecture is not nearly as satisfactory as the 
result for $\GLn$ (primarily because it's not necessary,
only claimed to be sufficient), but the situation for other groups
seems to be inherently less clean.

\setcounter{section}{0}

\section{Appendix: The equivalence of $GL_n$-honeycombs with a
definition of Berenstein-Zelevinsky patterns}\label{BZappendix}

Let $\eta$ be a hexagon with $120^\circ$ angles and two vertical edges.
Define the {\dfn torsion} of $\eta$ to be the length of the left edge
minus that of the right edge.

\begin{Proposition}
The torsion of $\eta$ and that of each $120^\circ$ rotation of $\eta$ agree.
\end{Proposition}

\begin{proof}
For a regular hexagon they are all zero. If one translates one edge of
$\eta$ out from the center, the edge shrinks and its two neighboring edges
grow, keeping the torsions equal. Any position of the hexagon can be
achieved by composing such translations.
\end{proof}

Let $h$ be a $GL_n$-honeycomb.  
We will assign a number to each region in the
$GL_n$ honeycomb tinkertoy $\tau_n$, 
other than the sectors at the three corners, using $h$. 
This will turn out to be a Berenstein-Zelevinsky pattern.

\begin{enumerate}
\item
 Each hexagon is assigned its torsion. 

\item
 Each \semiinfinite wedge on the NE long edge
of the honeycomb is assigned the length of its west edge.

\item
 Each \semiinfinite wedge on the NW long edge
of the honeycomb is assigned the length of its SE edge.

\item
 Each \semiinfinite wedge on the bottom long edge 
of the honeycomb is assigned the length of its NE edge.
\end{enumerate}

(This is set up so as to be $120^\circ$-rotation invariant.)

For any region not on the NW long edge,
the sum of the region-entries at and to the right of that point 
telescopes to the length of an edge, necessarily nonnegative. (Likewise
for $120^\circ$ rotations.)

\begin{figure}[Hht]
  \begin{center}
    \leavevmode
    \epsfig{file=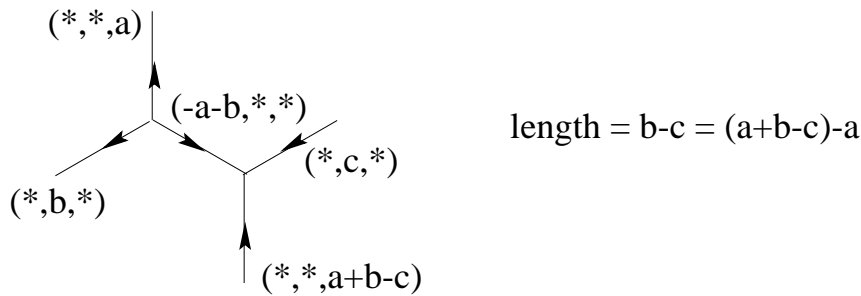,height=1.5in}
    \caption{Formula for the length of an edge. The edges in the 
      honeycomb are labeled with their constant coordinates.}
    \label{edgelength}
  \end{center}
\end{figure}

To determine the sum across an entire row is a little trickier.  We
need to relate the length of an edge to the constant coordinates on
neighboring edges. Rotate the edge
to align it with the finite edge in figure \ref{edgelength}. 

So the sum across an entire left-right row, with a \semiinfinite wedge
at the left end, is 
\begin{enumerate}
\item the sum of the lengths of the two finite edges of that wedge
\item the difference of the constant coordinates of the \semiinfinite edges
of that wedge.
\end{enumerate}
(By $120^\circ$-rotational symmetry the same is true for sums in
other directions.)

In particular, if the constant coordinates of the \semiinfinite edges
are interpreted as the coefficients $\lambda_i,\mu_i,\nu_i$ of three
dominant weights (in nonincreasing order) as in the rest of the paper,
the labeling of the regions exactly matches the definition of
{\dfn Berenstein-Zelevinsky pattern} as given in \cite{Z}
(only one of many realizations, others to be found in the original
\cite{BZ}).

The central theorem in \cite{BZ} gives a formula for Littlewood-Richardson
coefficients as the number of BZ patterns with given boundary values.
Their formulation is more suited to $SL_n$ than $GL_n$ calculations, and
as such, they need to include a caveat that BZ patterns count the LR
coefficient only if the sum of the three weights is in the root lattice
of $SL_n$. (If the sum is not in the root lattice, 
the LR coefficient is obviously zero, 
but there are still likely to be many BZ patterns which now have no known
representation-theoretic meaning.)
In the $GL_n$-adapted formulation of this paper this caveat 
does not appear.

For us, the BZ theorem reads as follows:

\begin{Theorem}\label{BZequiv}
Let $\lambda,\mu,\nu$ be a triple of dominant weights of $\GLn$.
Then the number of lattice $GL_n$-honeycombs 
whose \semiinfinite edges, indexed clockwise from the southwest,
have constant coordinates
$\lambda_1,\ldots,\lambda_n,\mu_1,\ldots,\mu_n,\nu_1,\ldots,\nu_n$
is the corresponding Littlewood-Richardson coefficient
$\dim (V_\lambda \tensor V_\mu \tensor V_\nu)^{\GLn}$.
\end{Theorem}

\section{Appendix: The hive model}\label{hiveAppendix}

In this section we introduce another model of the points in the
Berenstein-Zelevinsky cone, much closer to the BZ models in feel, 
that has the honeycomb-like property of having only ``local'' inequalities. 
It is not strictly necessary for the logic of the paper, 
but is very useful as an alternate model;
this is particularly true if one wants to
actually count tensor product multiplicities, rather than merely
prove them positive. The paper \cite{Bu} exposing our work
takes this model as the fundamental one.

As was mentioned elsewhere, the primary advantage of the honeycomb model
over the original BZ models is the naturality of the ``overlay'' operation;
in this paper this is only used in a sort of local way, when we
elide simple degeneracies. But this comes with a cost -- the linear
structure on the space of honeycombs is a bit difficult 
to see geometrically. In addition, the degrees of freedom of the honeycomb
are distinctly more obscure than in the BZ models.\footnote{%
In some sense, though, the strength of the honeycomb model as used in this 
paper that one can study the degrees of freedom left 
{\em while keeping some of the inequalities pressed,} 
which is not so easy to do in the BZ models.}
(And then there is the typesetting problem.)

Given $\tau$ a honeycomb tinkertoy, 
recall the dual graph $D(\tau)$ defined in section \ref{sec:degen}.
Define the vector space $\hive(\tau)$ to be labelings of
the vertices of $D(\tau)$ by real numbers.
This vector space naturally contains the lattice of integer labelings.

Recall that the embedded graph 
$D(\tau)$ is a collection of triangles, which we refer to as the 
{\dfn hive triangles}.
Of most interest to us are the rhombi formed by pairs of adjacent 
hive triangles. Each such rhombus has two {\dfn acute vertices} and
two {\dfn obtuse vertices}. There are three possible directions a
rhombus may face.

\begin{figure}[Hht]
  \begin{center}
    \leavevmode
    \epsfig{file=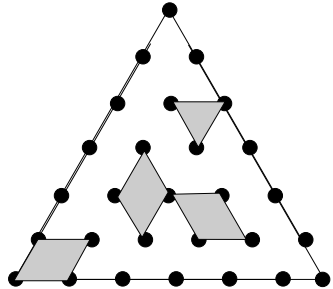,height=1.5in}
    \caption{The dual graph $D(\tau)$ corresponding to the
      $GL_6$ honeycomb tinkertoy, a little
      hive triangle, and a rhombus in each orientation.}
    \label{fig:6hive}
  \end{center}
\end{figure}

\begin{figure}[Hhtb]
  \begin{center}
    \leavevmode
    \epsfig{file=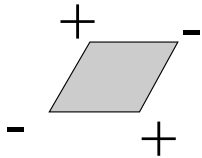,height=.8in}
    \caption{The dot-product picture of a rhombus inequality.}
    \label{fig:rhombussigns}
  \end{center}
\end{figure}

Each rhombus $\rho \subseteq H$ gives a functional on $\hive(\tau)$, 
defined as the sum at the obtuse vertices minus the sum at 
the acute vertices (as seen in figure \ref{fig:rhombussigns}).
This gives a {\dfn rhombus inequality}, asking that the rhombus
functional be nonnegative. 
The cone in $\hive(\tau)$ satisfying all these inequalities we
denote $\HIVE(\tau)$, and we will call its elements $\tau$-hives. 

\begin{Proposition}
  Let $\tau$ be a honeycomb tinkertoy. Then there is a
  $\integers$-linear correspondence between configurations of $\tau$,
  and $\tau$-hives whose leftmost top entry is zero, in such a way that
  the constant coordinates of the boundary edges of the honeycomb are
  differences of boundary entries on the corresponding $\tau$-hive.
\end{Proposition}

\begin{proof}
Start with a configuration $h$ of $\tau$.  We label the vertices in
$D(\tau)$ inductively, starting with a zero in the leftmost top entry, 
and filling in
as follows: whenever we move southwest or east, we increase the value
by the constant coordinate of the edge crossed; southeast, we decrease
by that constant coordinate.

Our first worry is that different paths will cause us to try to fill
different numbers in the same hexagon. That this doesn't happen is a simple 
consequence of the sum-equals-zero property at a vertex of the honeycomb.

Second, we need to know that the result is a hive. Not surprisingly, the
rhombus inequalities are equivalent to the edge lengths being nonnegative.

Lastly, since we define the hive entries by inductively adding up
constant coordinates of edges, the boundary of the hive naturally ends
up being the partial sums of those constant coordinates.
(The sum-equals-zero property is involved in seeing this for 
some of the boundary edges.)
\end{proof}

Combining this with the theorem in the appendix relating $GL_n$-honeycombs
to Littlewood-Richardson coefficients, we find that if $\lambda,\mu,\nu$ 
are integral, the number of hives with boundary formed from partial
sums of $\lambda,\mu,\nu$ is a Littlewood-Richardson coefficient.
(In \cite{Bu}, there is given a simple bijection between hives and a
standard formulation of the Littlewood-Richardson rule.)

There is a pleasant geometric way to interpret the rhombus inequalities.
Extend the hive to a piecewise linear function, affine-linear on each
little hive triangle. Then the rhombus inequalities state that this 
function is convex. Each rhombus {\em equality} says that the function
is actually linear across the boundary down the middle of the rhombus,
i.e. that the regions on which the function is affine-linear are 
larger than just the little hive triangles. 

In this way, the set of tight
rhombus inequalities determines a certain decomposition of the
convex region in $\B$ bearing $D(\tau)$ into regions 
(which we dub ``flatspaces'' due to the geometric interpretation above). 
This is exactly the decomposition
into the regions of the degeneracy graph from section \ref{sec:degen}.
We mentioned there that the degeneracy graph remembers only the
``combinatorial information'' about a honeycomb;
we see now that it is the hive that finishes the job.

Again, we refer readers to the honeycomb/hive applet to see these hives
and convex graphs in action, at
\begin{center}
{\tt http://www.alumni.caltech.edu/\~{}allenk/java/honeycombs.html.}  
\end{center}

\bibliographystyle{alpha}

\begin{thebibliography}{10}

\bibitem[Bu]{Bu} A. Buch, The saturation conjecture (after A. Knutson and
T. Tao), notes from a talk at Berkeley September 1998.

\bibitem[BS]{BS} L. Billera, B. Sturmfels,
Fiber polytopes,
{\it Annals of Math.} {\bf  135} (1992), no. 3, 527--549. 

\bibitem[BZ]{BZ} A.~Berenstein, A.~Zelevinsky, Triple multiplicities for
$sl(r+1)$ and the spectrum of the exterior algebra of the adjoint
representation, {\it J. Alg. Comb.}, {\bf 1} (1992), 7 - 22. 

\bibitem[E]{E} A.G.~Elashvili, Invariant algebras, 
{\it Advances in Soviet Math.}, {\bf 8} (1992), 57-64. 

\bibitem[F]{F} W. Fulton,
Eigenvalues of sums of Hermitian matrices (after A. Klyachko), 
S\'eminaire Bourbaki. (1998). 

\bibitem[FH]{FH} W. Fulton, J. Harris,
Representation theory, Springer-Verlag (1991).

\bibitem[GP]{GP} O. Gleizer, A. Postnikov,
Littlewood-Richardson coefficients via Yang-Baxter equation,
in preparation.

\bibitem[GZ]{GZ} V. Guillemin, C. Zara,
Equivariant de Rham theory and graphs, {\tt math.DG/9808135}.

\bibitem[H]{Horn} A.~Horn, Eigenvalues of sums of Hermitian matrices,
{\it Pacific J. Math.}, {\bf 12} (1962), 225-241.

\bibitem[Hon2]{Hon2} A. Knutson, T. Tao, C. Woodward,
The honeycomb model of $\GLn$ tensor products II: 
Facets of the L-R cone, in preparation.

\bibitem[KSZ]{KSZ}
M. Kapranov, B. Sturmfels, A. Zelevinsky, 
Quotients of toric varieties. 
{\it Math. Annalen.} {\bf 290} (1991), no. 4, 643--655.

\bibitem[Kl]{Kly} A.A.~Klyachko, Stable vector bundles and Hermitian
operators, {\it IGM, University of Marne-la-Vallee preprint} (1994).

\bibitem[KMP]{KM}
S. Kumar,
Proof of the Parthasarathy-Ranga Rao-Varadarajan conjecture. 
{\it Invent. math.} {\bf 93} (1988), no. 1, 117--130. 

O. Mathieu,
Construction d'un groupe de Kac-Moody et applications.
{\it Compositio Math.} {\bf 69} (1989), no. 1, 37--60. 

P. Polo,
Vari\'et\'es de Schubert et filtrations excellentes.
{\it Ast\'erisque.} {\bf 10-11} (1989) 281-311.

\bibitem[Ku]{Ku} G. Kuperberg,
Spiders for rank two Lie algebras,
{\tt math.QA/9712143},
{\it Comm. Math. Phys.}, 180(1):109-151, 1996

\bibitem[MFK, chapter 8]{MFK}
D. Mumford, J. Fogarty, F. Kirwan,
Geometric invariant theory. Third edition.
{\it Ergebnisse der Mathematik und ihrer Grenzgebiete},
Springer-Verlag, 1994.

\bibitem[Ze]{Z}  A.~Zelevinsky, Littlewood-Richardson semigroups,
{\tt math.CO/9704228}.

\bibitem[Zi, pp. 293-4]{Zi}
G. Ziegler,
Lectures on polytopes,
{\it Graduate Texts in Mathematics}, 152. 
Springer-Verlag, 1995.

\end{thebibliography}

\end{document}